\documentclass[aos,MSNbibl,nameyear,dvips]{arximspdf}

\doi{10.1214/13-AOS1085} 
\volume{41}
\issue{2}
\pubyear{2013}
\firstpage{536}
\lastpage{567}

\makeatletter

\newcommand{\cal}{\mathcal}

\newtheorem{theorem}{Theorem}[section]
\newtheorem{lemma}{Lemma}[section]

\newproclaim{remark}{Remark}[section]
\newproclaim{condition}{Condition}[section]

\newcommand{\R}{\mathbb{R}}

\newcommand{\Var}{\operatorname{Var}}

\newcommand{\PP}{\mathbb{P}}

\newcommand{\EE}{\mathbb{E}}
\newcommand{\pa}{\operatorname{pa}}

\makeatother

\begin{document}
\begin{frontmatter}

\title{$\ell_0$-penalized maximum likelihood for
sparse directed acyclic graphs}
\runtitle{Penalized maximum likelihood for DAG}

\begin{aug}
\author[A]{\fnms{Sara} \snm{van de Geer}\corref{}\ead[label=sg]{geer@stat.math.ethz.ch}}
\and
\author[A]{\fnms{Peter} \snm{B\"uhlmann} \ead[label=pb]{buhlmann@stat.math.ethz.ch}}
\runauthor{S. van de Geer and P. B\"uhlmann}
\affiliation{ETH Z\"urich}
\address[A]{Seminar f\"ur Statistik\\
ETH Z\"urich\\
Switzerland\\
\printead{sg}\\
\hphantom{E-mail: }\printead*{pb}} 
\end{aug}

\received{\smonth{5} \syear{2012}}
\revised{\smonth{11} \syear{2012}}

%
\begin{abstract}
We consider the problem of regularized maximum likelihood estimation
for the structure and parameters of a high-dimensional, sparse directed
acyclic graphical (DAG) model with Gaussian distribution, or
equivalently, of a Gaussian structural equation model. We show that the
$\ell_0$-penalized maximum likelihood estimator of a DAG has about the
same number of edges as the minimal-edge I-MAP (a DAG with minimal
number of edges representing the distribution), and that it converges
in Frobenius norm. We allow the number of nodes $p$ to be much larger
than sample size $n$ but assume a sparsity condition and that any
representation of the true DAG has at least a fixed proportion of its
nonzero edge weights above the noise level. Our results do not rely on
the faithfulness assumption nor on the restrictive strong faithfulness
condition which are required for methods based on conditional
independence testing such as the PC-algorithm.
\end{abstract}

%

\end{frontmatter}

\section{Introduction}

Directed acyclic graphs (DAGs) and corresponding directed graphical models
are key concepts for causal inference; see, for example, the books by
\citet{sgs00} and \citet{pearl00}. From an estimation point of view, a
first step consists of estimating the Markov equivalence class of the true
underlying causal DAG based on observational data, and from this, one
can infer identifiable causal effects and lower bounds for all causal
effects [\citet{makapb09}]. This strategy has been applied to, and to a certain
extent validated using high-throughput, and hence high-dimensional,
data in
biology
[\citet{mapb10}]. It is of primary importance to understand
limitations and potential of methods in terms of subtle and often
uncheckable assumptions, and in this respect, our results here shed new
light.

We focus here on the problem of estimating the
Markov equivalence class of DAGs, or more generally of a so-called
minimal-edge I-MAP, in the setting of observational Gaussian data where the
number $p$ of variables or nodes in the DAG may greatly
exceed sample size $n$. We consider the $\ell_0$-penalized maximum likelihood
estimator, and we relate and compare our new results and conditions to the
popular PC-algorithm [\citet{sgs00}] and its theoretical analysis. To the
best of our knowledge, the latter is so far the
only work providing theoretical guarantees for inferring the Markov equivalence
class of DAGs in the high-dimensional setting. We
emphasize that the popular $\ell_1$-norm regularization of the likelihood
is inappropriate here, leading to an objective function which is not
constant over
equivalent DAGs encoding the same distribution. On
the other hand, $\ell_0$-penalization leads to invariant scores over
equivalent DAGs. The computational difficulties are primarily due to the
nonconvex constraint that the directed graph should be acyclic, and the
additional issue with the $\ell_0$- in comparison to, for example, $\ell_1$-norm
penalization is, surprisingly, rather in favor of the former. A~computationally feasible algorithm for
exact $\ell_0$-penalized maximum
likelihood estimation for the Markov equivalence class of DAGs has been
proposed by \citet{sil06}; for larger graphs, with more than about 20 nodes,
approximate algorithms can be used [\citet{chick02,hapb11}]. More
details are given in Section~\ref{subsecell0-pen}.

\subsection{Relation to other work}

We analyze the $\ell_0$-penalized maximum likelihood estimator for the
equivalence class of DAGs in the Gaussian setting. Pioneering work for the
low-dimensional case on this problem has been done by \citet{chick02} who
proved consistency of the BIC-score and provided an algorithm, called
greedy equivalent search (GES), which greedily proceeds in the space of
Markov equivalence classes. While the GES-algorithm can also be used in the
high-dimensional scenario [\citet{hapb11}], the asymptotic consistency
of BIC
is established only for the case with a fixed distribution (with $p <
\infty$) where the
sample size $n \to\infty$. Chickering's first significant analysis does
not provide any insights for the high-dimensional case with its subtle
interplay of signal strength, noise level and identifiability conditions.

\citet{robins03} present refined analyses for causal inference under the
view point of uniform consistency as sample size $n \to\infty$. There,
problematic issues with the so-called faithfulness condition (see Section
\ref{subsecDAGSEM}) arise, and \citet{zhang03} introduce the notion of
strong faithfulness [see (\ref{strong-faith})], as a way to address
some of
the raised major problems. None of
these works consider high-dimensional inference, but their pointing to the
faithfulness condition and its version are important.

\citet{kabu07} provide consistency results of the
PC-algorithm [\citet{sgs00}] for estimating the Markov equivalence class of
DAGs based on Gaussian observational data, in the high-dimensional, sparse
setting. One of the conditions used is a restricted version of strong
faithfulness in (\ref{strong-faith}). Our analysis with the penalized
MLE is completely different and circumvents faithfulness and strong faithfulness
conditions which are often very restrictive as shown by
\citet{uhleretal12}; see also below.

The theoretical high-dimensional analysis presented here is very different
and more challenging than for multivariate regression or covariance
estimation, due to the \textit{unknown} order among the variables. For known
order, as, for example, in time series problem, \citet{shojaie10}
present results
for estimation of high-dimensional DAGs; but the case with unknown order
considered in the present paper requires major new theoretical ideas and
development.

\subsection{Directed graphical and structural equation
models}\label{subsecDAGSEM}

Consider the following model. There is a DAG $D_0$ whose $p$ nodes
correspond to random variables $X_1,\ldots,X_p$: assume that
%
\begin{eqnarray}
\label{modDAG}
&\displaystyle
X_1,\ldots,X_p \sim{\cal
N}_p(0,\Sigma_0) \mbox{ with density } f_{\Sigma_0}(
\cdot),&
\nonumber\\[-8pt]\\[-8pt]
&\displaystyle {\cal N}_p(0,\Sigma_0) \mbox{ is Markovian with
respect to $D_0$},&\nonumber
\end{eqnarray}
where the Markov property can be understood as the factorization property
where the joint Gaussian density $f_{\Sigma_0}(x_1,\ldots,x_p)$ can be
factorized as
\[
f_{\Sigma_0}(x_1,\ldots,x_p) = \prod
_{j=1}^p f_{\Sigma_0}(x_j|x_{\pa(j)})
\]
with $\pa(j)$ denoting the set of parents of node $j$; cf. \citet{lauritzen96}.

It is well known that, in general, there exists another DAG $D$ such
that the
distribution ${\cal N}(0,\Sigma_0)$ is Markovian with respect to
$D$. Assuming faithfulness (see below), the set of all such other DAGs
build the so-called Markov equivalence class
${\cal E}(D_0)$ which can be characterized in terms of a chain graph with
undirected and directed edges; cf. \citet{ander97}. The Markov equivalence
class ${\cal E}(D_0)$ can be identified from the observational data
distribution ${\cal N}(0,\Sigma_0)$ under the assumption of
faithfulness.\vspace*{9pt}

\textit{Definition of faithfulness}; \textit{cf. Spirtes}, \textit{Glymour and
Scheines}
(\citeyear{sgs00}): For a DAG~$D$, a
distribution $P$ is
called \textit{faithful} with respect to $D$ if and only if all conditional
independences are encoded by the DAG $D$.\vspace*{9pt}

Faithfulness is stronger than a Markov assumption: the latter allows us to
infer some conditional independences from the DAG while the former requires
that \textit{all} of them can be inferred from the DAG (i.e., only the ones
which are entailed by a Markov condition). Failure of faithfulness is
``rare,'' having Lebesgue measure zero, if the edge weights (the
coefficients in the equivalent linear structural equation model) are chosen
from a distribution which is absolutely continuous with respect to Lebesgue
measure.\vadjust{\goodbreak} However, for statistical estimation, we often require
sufficiently strong detectability of conditional dependencies, given by the
notion of strong faithfulness.\vspace*{9pt}

\textit{Definition of strong faithfulness in the Gaussian case} [\textit{Zhang and
Spirtes} (\citeyear{zhang03})]:
For a DAG $D$, a
Gaussian distribution $P$ is called $\tau$-\textit{strongly faithful} with
respect to $D$ if and only if
%
\begin{eqnarray}
\label{strong-faith}
&&P \mbox{ is faithful with respect to $D$}\quad\mbox{and}
\nonumber
\\
& &\min\bigl\{\bigl|\mathrm{Parcorr}(X_j,X_k|X_S)\bigr|;
\operatorname{Parcorr}(X_j,X_k|X_S) \neq0,
\\
& &\hspace*{24.4pt} S \subseteq\{1,\ldots,p\}\setminus\{j,k\}, j,k \in\{1,\ldots, p\}\
(j \neq
k)\bigr\} \ge\tau.\nonumber
\end{eqnarray}

A typical requirement is strong faithfulness with $\tau\asymp
\sqrt{\mbox{sparsity} \cdot\log(p)/n}$, see also below for the
PC-algorithm. Strong
faithfulness can be
viewed as a condition of ``signal strength'' in terms of nonzero partial
correlations. As shown in \citet{uhleretal12}, strong faithfulness is a
very restrictive condition for many DAGs, and the same applies to a
slightly weaker restricted strong faithfulness assumption; cf.
\citet{uhleretal12}; see also Section
\ref{subsecbetamin-strongfaith}. At the same time, it is
essentially unavoidable for any algorithm for inferring the Markov
equivalence class ${\cal E}(D_0)$ which relies on conditional independence
testing. The most prominent example is the PC-algorithm [\citet{sgs00}]:
consistent estimation for the Markov equivalence class ${\cal E}(D_0)$ is
proved in \citet{kabu07} for the Gaussian case assuming a strong
faithfulness condition.
The results in this paper for
the $\ell_0$-penalized MLE do not require a strong faithfulness condition
as in (\ref{strong-faith}) (nor the slightly weaker restricted strong
faithfulness condition): the reason is that
the method is \textit{not} relying on conditional independence testing but
rather on penalized parameter estimation in terms of a linear structural
equation model, as we explain next.

A Gaussian DAG model as in (\ref{modDAG}) can always be equivalently
represented as a linear structural equation model,
%
\begin{equation}
\label{modSEM} X_j = \sum_{k \in\pa(j)}
\beta^0_{kj} X_k + \epsilon_j\qquad
(j=1,\ldots,p),
\end{equation}
where $\epsilon_1,\ldots, \epsilon_n$ are independent,
$\epsilon_j \sim{\cal N}(0,|\omega_j^0| ^2)$ and $\epsilon_j$
independent of
$\{X_k; k \in\pa{(j)}\}$; note that $\pa(j) = \pa_{D_0}(j)$ depends
on the
true DAG $D_0$.

\section{The setting and the estimator}

We use here and in the sequel a terminology which does not rely on the
standard language from graphical modeling since the required basics for the
Gaussian case [see models (\ref{modDAG}) and (\ref{modSEM})] can be
developed in a straightforward way.

We consider $n$ i.i.d. observations from the structural equation model
(\ref{modSEM}) which is equivalent to model (\ref{modDAG}).
We denote by $X:= (X_1,\ldots, X_p)$ the $n \times p $ data matrix with
$n$ i.i.d. rows, each of them being ${\cal N} ( 0, \Sigma_0
)$-distributed, where $\Sigma_0$ is a nonsingular covariance matrix.
The relations between the
variables in a row can be represented as
%
\begin{equation}
\label{model} X = X B_0 + E,
\end{equation}
where $B_0:= ( \beta_{k,j}^0 ) $ is a $p \times p $ matrix with
$\beta_{j,j}^0 = 0 $ for all $j$, and where
$ E $ as an $n \times p $ matrix of noise vectors $E:= ( \epsilon_1,\ldots, \epsilon_p)$, with $\epsilon_j$ independent of $X_k$ whenever
$\beta_{k,j}^0 \not= 0 $.\footnote{Note\vspace*{1pt} that in relation to the true DAG
$D_0$ in model (\ref{modSEM}), $\beta_{k,j}^0 = 0$ for $k \notin
\pa_{D_0}(j)$. We do not make such explicit constraints here since we aim
for a smallest DAG representing the distribution of $X$.}
Furthermore, $E$ has $n$ i.i.d. rows which are ${\cal N} ( 0,
\Omega_0)$-distributed, with
$\Omega_0:= \operatorname{diag} ( | \omega_1^0 |^2,\ldots, |\omega_p^0 |^2
)$ a
$p \times p $ diagonal matrix.

Model (\ref{model}) implies that
\[
\Sigma_0 = \bigl[ ( I-B_0)^{-1}
\bigr]^T \Omega_0 \bigl[( I- B_0)^{-1}
\bigr].
\]
We call $(B_0, \Omega_0)$ a DAG corresponding to
$\Sigma_0$.\footnote{This deviates from the classical definition where
the DAG is only a (directed acyclic) graph; we use the short
terminology ``DAG'' for the whole graphical model with the distribution
and the graph encoded by the coefficient matrix $B$ and the error
variances $\Omega$.} The set of edges of this DAG is denoted by $s_0:=
s_{B_0}:= \{ (k,j)\dvtx  \beta_{k,j}^0 \not= 0 \} $, and in
fact,\vspace*{1pt}
$\beta_{k,j}^0 \neq0$ encodes for a directed edge $k \to j$. We will
assume in Condition~\ref{edgescondition} that $(B_0,\Omega_0)$ is
sparse, in the sense that its number of (directed) edges $s_0$ is
small.

As we described in Section~\ref{subsecDAGSEM} with the concept of a Markov
equivalence class, there are several DAGs $( \tilde B_0, \tilde\Omega_0)$
describing the same $\Sigma_0$ and thus the same Gaussian distribution $P=
{\cal N} (0, \Sigma_0)$. Throughout this paper, $(B_0, \Omega_0)$ is defined
as a DAG with a minimal number of edges (and it may not be unique). We call
such a DAG a minimal-edge
I-MAP.\footnote{It is a minimal I-MAP [\citet{sgs00}, Section 2.3.1]
with the
additional property that it has minimal number of edges.}
%
\begin{remark}\label{equivalentremark}
We call
two DAGs, $(B_1, \Omega_1)$ and $(B_2, \Omega_2)$, equivalent if
$\Theta(B_1, \Omega_1)= \Theta( B_2, \Omega_1)$ and
if in addition they have the same number of edges.\footnote{This
definition is not the same as for a Markov equivalence class. Assuming
faithfulness, both definitions coincide.}
In our analysis, we will identify DAGs which are in the same
equivalence class.
Thus, our aim is to estimate an arbitrary member of the equivalence class
of $(B_0, \Omega_0)$, by a suitable member in the equivalence class of
$(\hat B, \hat\Omega_0)$.
\end{remark}

\subsection{\texorpdfstring{The $\ell_0$-penalized maximum likelihood estimator}{The l0-penalized maximum likelihood estimator}}\label{subsecell0-pen}
We use a penalized maximum likelihood
procedure to estimate the DAG $(B_0, \Omega_0)$. Let
\[
\Sigma_n:= X^T X / n
\]
be the empirical covariance matrix based on the observations $X$.
Given a $p \times p$ nonsingular covariance matrix $\Sigma$,\vadjust{\goodbreak} with inverse
$\Theta:= \Sigma^{-1}$, the minus log-likelihood is proportional to
\[
l_n (\Theta):= \operatorname{trace} ( \Theta\Sigma_n ) - \log
\operatorname{det} ( \Theta).
\]
We consider inverse covariance matrices that can be represented
as a DAG. That is, we let
\[
\Theta:= \Theta( B, \Omega):= (I-B) \Omega^{-1} (I-B)^T,
\]
where $( B, \Omega)$ is a DAG. The latter means that $\Omega$ is a
positive diagonal matrix and that, up to a permutation $\pi$ of
the rows and columns, $B$ can be written as a lower-diagonal matrix.

The $\ell_0$-penalized maximum likelihood estimator is
%
\begin{eqnarray}
\label{ell0MLE} \hat\Theta &=& \mathop{\arg\min}_{B, \Omega} \bigl\{l_n
( \Theta) + \lambda^2 s_{B}\dvtx
\Theta = \Theta( B, \Omega),\nonumber\\[-8pt]\\[-8pt]
&&\hspace*{38.3pt}\mbox{for some DAG } ( B, \Omega)
\mbox{ with } B
\in{\cal B} \bigr\}.\nonumber
\end{eqnarray}
Here $s_{B}$ is the number of nonzero elements in $B$ (corresponding to
the number of edges in the DAG) and $\lambda\ge0$ is a tuning
parameter.
The estimator is denoted by
$\hat\Theta:= \Theta( \hat B, \hat\Omega) $.
It has $\hat s:= s_{\hat B} $
edges.
The collection ${\cal B}$ is the set of all edge weights $B$ of DAGs
$(B, \Omega)$
which have at most
$\alpha n/ \log p $ incoming edges (parents) at each node, where $\alpha>0$
is given (see Condition~\ref{hatedgescondition}), or a subset thereof.
We will discuss this restriction in Section~\ref{edgessection}.
We throughout assume
$B_0 \in{\cal B}$, that is, that the restrictions one puts on the edge
weights are correct.

The $\ell_0$-penalty in the estimator ensures that the penalized likelihood
remains the same among all equivalent representations, for example,
among all DAGs
from the same Markov equivalence class or among the equivalence class
described in Remark~\ref{equivalentremark} above. This would not be true
when choosing, for example, an $\ell_1$-norm penalization $\|B\|_1:=
\sum_{k,j} |\beta_{k,j}|$.
From a computational point of view, the
main difficulty is the optimization over the space of DAGs ${\cal B}$:
current algorithms are tailored for the $\ell_0$-penalty (see next
paragraph), and in
this sense, optimization of the $\ell_0$-penalized log-likelihood is better
tractable than using an $\ell_1$-norm regularization. For
problems with up to about $p \approx20$ nodes, dynamic programming can be
used [\citet{sil06}], while for large-scale applications, greedy equivalent
search has been reported to perform well [\citet{chick02,hapb11}].

The dynamic programming method [\citet{sil06}] is optimizing the
$\ell_0$-penalized negative log-likelihood in (\ref{ell0MLE}) over the
space of all DAGs with $B \in{\cal B}$. It crucially relies on the
assumption that the objective function, called the score, can be decomposed
locally such that $\operatorname{score}(D) = \sum_{j=1}^p
g(X_j,X_{\mathrm{pa}_D(j)})$ for some function $g(\cdot)$, where
$g(X_j,X_{\mathrm{pa}_D(j)})$ involves only data from the variables
$j$ and
$\mathrm{pa}_D(j)$. The $\ell_0$-penalized score is decomposable, whereas,\vadjust{\goodbreak}
say, $\ell_1$-norm penalization does not lead to a decomposable score. The
greedy equivalent search algorithms [\citet{chick02,hapb11}] are greedily
optimizing the objective function in (\ref{ell0MLE}) in the space of
equivalence classes, a much smaller space than the space of DAGs. The
greedy steps are forward, backward and turning edges moves such that the
score is improved most when stepping from one Markov equivalence class to
another one: this can be done very efficiently without enumerating all DAG
members in the equivalence class. Such greedy equivalent search
algorithms rely crucially on the fact that the objective score
function is
constant within an equivalence class and that the score is decomposable:
again, this is true with $\ell_0$-penalization but fails for
$\ell_1$-norm regularization.

The previous paragraph already gives a good reason why
$\ell_0$-penalization is to be preferred over the $\ell_1$-norm analogue,
namely, that the computational techniques are tailored and
surprisingly easier for the former than the latter ($\ell_1$-norm
regularization would require to search over the whole space of DAGs, and
the dynamic programming trick mentioned above cannot be used). Furthermore,
$\ell_1$-norm regularization is usually motivated as a convex relaxation,
and this is not true in our context since the DAG constraint in
(\ref{ell0MLE}) for the matrix $B \in{\cal B}$ remains to
cause that the optimization problem is nonconvex. In addition, the price
to pay, in a standard regression problem with an $\ell_1$-norm
regularization, is a bias which should be further controlled with, for
example, adaptive $\ell_1$-norm regularization [\citet{zou06}] or thresholding
[\citet{Geer11}]. We have not considered a theoretical analysis of the
$\ell_1$-norm penalized maximum likelihood estimator
of a DAG. We believe this not to be more difficult than the $\ell_0$-norm
penalized maximum likelihood estimator although the proofs (see
Section~\ref{proofsoutlinesection} for an outline) may need
different arguments.

\subsubsection{The main results and their implications}\label{subsecmainresults}
We show in Theorem~\ref{maintheorem} that with a choice of the tuning
parameter $\lambda^2$ of order $\log p / n (( p/ s_0) \vee1 )$, the
number of edges $\hat s$ of the estimator is of the same
order of magnitude as the number of edges $s_0$ of a true underlying
minimal-edge I-MAP. Moreover, we show that $\hat B$ and $\hat\Omega$
converges in Frobenius norm to some $\tilde B_0$ and $\tilde\Omega_0$,
respectively, where
$\Theta(\tilde B_0, \tilde\Omega_0) = \Theta_0$ is a representation of
the true DAG with $\tilde s$ edges (see Section~\ref{subsecpermutations}),
and with $\tilde s$ of the same order of magnitude as $s_0$.
The rate of convergence is of order~$\lambda^2 s_0 $. To arrive at this result,
we need that at least a fixed proportion of the nonzero coefficients
of any
representation of the true DAG is above the ``noise level,'' the latter
being of order $\sqrt{\log p / n}(\sqrt{p/s_0} \vee1) $ (see Condition
\ref{beta-mincondition}): in analogy to regression, we call this the
``beta-min'' condition. Some of our other conditions are trivially
satisfied when $p = {\mathcal O} (n/\log(n))$ is sufficiently small.

The ``noise level'' indicates two regimes for $(n,p,s_0)$. If $s_0$ is at
least of the order of $p$ (or larger), then the ``noise level'' is of the
order $\sqrt{\log(p)/n}$ which is small even if $p$ is very large relative
to $n$. This scenario is often realistic saying that a fixed nonzero proportion
of the nodes
has at least one parent: we call it the standard sparsity regime. The other
scenario is for $s_0 \ll p$, corresponding to very sparse DAGs, which we
call the ultra-high sparsity regime. The reason for the two different
noise level scenarios is that we estimate $p$ error variances in
$\Omega_0$: when $s_0 \ll p$, the term from estimating these error
variances is dominating.

When making a more stringent ``beta-min'' condition and choosing the
regularization parameter of larger order than the ``noise level'' [or the
same order if $s_0 = {\mathcal O} (1)$], we can recover the minimal-edge
I-MAP.

The $\ell_0$-penalized MLE can be easily adapted to the case where the
noise variances in $\Omega_0$ are known up to a scalar, for example, when
all noise variances are known to be equal but their value is unknown. Then,
the true DAG can be identified from the distribution
[\citet{Peters2012}]. We show that in this case, under an identifiability
condition on the noise variances, the $\ell_0$-penalized maximum likelihood
estimator finds the representation (and hence the true DAG) with the
prescribed noise variances, and the rate of convergence for the Frobenius
norm is of order $\lambda^2 s_0$; see Theorem~\ref{equalvariancestheorem}.
We assume in this context that $p$ is sufficiently smaller than $n/ \log
n$.

\begin{remark}\label{remark-markoveqclass}
If the minimal-edge I-MAP can be inferred and assuming in addition the
faithfulness condition (but \textit{not}
requiring strong faithfulness) (see Section~\ref{subsecDAGSEM}), we can
recover the true underlying Markov equivalence class. This then allows
us to
derive bounds for causal inference, exactly along the lines of
\citet{makapb09}.
\end{remark}

\subsection{Permutations and the order of the variables}\label{subsecpermutations}
Model (\ref{model}) can be written as
\[
X_j = X \beta_j^0 + \epsilon_j,\qquad
j=1,\ldots, p,
\]
with $\beta_j^0 $ the $j$th column of $B_0$.
Let us write for any vector $\beta\in\R^p $,
\[
\| X \beta\|^2:= \beta^T \Sigma_0 \beta,\qquad \|
X \beta\|_n^2:= \beta^T \Sigma_n
\beta.
\]
For a permutation $\pi$ of $\{ 1,\ldots, p \}$, which plays the role of
an order of the variables, we let $\tilde B_0
(\pi)$ be the
matrix obtained by doing a Gram--Schmidt orthogonalization for \mbox{$\|
\cdot\|$},
starting with $X_{\pi_p}$ and finishing\vspace*{-1pt} by projecting $X_{\pi_1}$ on
$X_{\pi_2}, \ldots, X_{\pi_p}$. Moreover, we let $\tilde\Omega_0 (\pi
)$ be the
diagonal matrix of the error variances. Note thus that $\tilde B_0 (\pi
)$ is
lower-diagonal after permutation of its rows and columns. Furthermore,
\[
\Theta_0 = \Theta_0 \bigl( \tilde B_0 (
\pi), \tilde\Omega_0 ( \pi)\bigr)\qquad \forall\pi.
\]

The set\vspace*{2pt} of incoming edges at node $j$ [nonzero coefficients of the $j$th
column of $\tilde B_0(\pi)$] is denoted
by $\tilde S_j ( \pi)$, and we let $\tilde s_j (\pi):= | \tilde S_j (\pi)|$.\vadjust{\goodbreak}
Moreover, we let $\tilde s( \pi)= \sum_{j=1}^p \tilde s_j (\pi)$ be the total
number of edges of $\tilde B_0 (\pi)$. Thus, \mbox{$\tilde s (\pi) =
s_{\tilde B_0 (\pi) } $}.

Let $\hat B$---or one of the members in its equivalence class---be
lower-diagonal after permutation $\hat\pi$, and define
$\tilde B_0:= B_0 ( \hat\pi)$.
The number of edges of $\tilde B_0 $ is denoted by
$\tilde s = \tilde s ( \hat\pi)$.
The DAGs $(\hat B, \hat\Omega)$ and $( \tilde B_0, \tilde\Omega_0)$
share the same lower-diagonal structure (but not necessarily the same set
of zero coefficients). We will show that $\tilde s:= \tilde s ( \hat
\pi)$
is with large probability of the same order of magnitude as $s_0$; see
Theorem~\ref{maintheorem}. Thus if
the true DAG $(B_0, \Omega_0)$ is sparse, then with large probability the
DAG $( \tilde B_0 ( \hat\pi), \tilde\Omega_0 ( \hat\pi) )$ is sparse
as well, which means that on average, the number of incoming edges at a
node is small.

Note that $\hat\pi$ is a random permutation and that there are in
total a
large number (namely $p!$) permutations. Analytical control over these
$p!$ permutations requires a very different technique than dealing with
known order [\citet{shojaie10}] or with multivariate regression or covariance
estimation problems. We explain this in more detail in Section
\ref{T1section}.

\subsubsection{AR(1)-model as an example}\label{subsecAR1}
Suppose the true DAG is a directed chain from
$X_p$ along $X_{p-1},\ldots, X_2$ to $X_1$ with a corresponding
structural equation model [AR(1)-model],
\begin{eqnarray*}
X_p &=& \epsilon_p,
\\[-1pt]
X_j &=& \beta^0 X_{j+1} +
\epsilon_j\qquad (j=1,\ldots,p-1),
\end{eqnarray*}
where $\epsilon\sim{\cal N}(0,\Omega_0)$ with $\Omega_0 = \operatorname{diag}
(1 -
(\beta^0)^2,\ldots,1 - (\beta^0)^2, 1)$ and $|\beta^0| < 1$. The error
variances are chosen such that $\Var(X_j) = 1$ for all $j$. The covariance
matrix is of Toeplitz form $(\Sigma_0)_{ij} = (\beta^0)^{|i-j|}$ and
the model
satisfies the directed global Markov property which is equivalent to the
concept of d-separation; cf. \citet{lauritzen96}, Section 3.2.2.

Therefore, we have that projecting
\[
X_{\pi_j} \mbox{ on } X_{\pi_{j+1}},\ldots,X_{\pi_p}\qquad (j =1,\ldots, p-1)
\]
leads to at most two nonzero regression coefficients in every column of
$\tilde B_0 (\pi)$ (corresponding to the largest index $k_1 < \pi_j$ and
smallest index $k_2 > \pi_j$ if $\pi_{j+1},\ldots,\pi_p$ contains indices
smaller and larger than $\pi_j$; or corresponding to the largest $k <
\pi_j$ if $\pi_{j+1},\ldots,\pi_p$ contains only smaller indices than
$\pi_j$; or corresponding to the smallest $k > \pi_j$ if $\pi
_{j+1},\ldots,\pi_p$ contains only larger indices than $\pi_j$). Thus,
we have
that $\tilde{s}_j(\pi) \le2$ for all $j$ and all $\pi$, and hence Condition
\ref{edgescondition}, given below, holds.

The absolute\vspace*{1pt} values of the nonzero coefficients in $\tilde{B}_0(\pi) =
(\tilde{\beta}^0_{k,j}(\pi))$ decrease monotonely as the
index-distance $d(j) = \min_{k = j+1,\ldots,p} |\pi_k - \pi_j |$
increases. Thus,
for fixed $j$ and whenever $d(j) >
\Delta$ for some (large)\vspace*{1pt} value of~$\Delta$, there are at most two [since
$\tilde{s}_j(\pi) \le2$] coefficients with $|\tilde{\beta}^0_{k,j}(\pi)|
\le C(\Delta)$ for some value $C(\Delta)$\vadjust{\goodbreak} (which decreases as $\Delta$
increases and which depends on~$\beta^0$). Therefore, clearly, there
are at
most $2( \lfloor p /\Delta
\rfloor+ 1)$ coefficients (edges)\vspace*{1pt} whose values satisfy
$|\tilde{\beta}^0_{k,j}(\pi)| \le C(\Delta)$, and all\vspace*{1pt} other nonzero
coefficients (at least $p-\lfloor p/\Delta\rfloor-2$)\footnote{There are at least $p -(\lfloor p/\Delta\rfloor+1)$ indices
(nodes) $j$ with $d(j) \le\Delta$; and there is at least one nonzero
coefficient (edge) from all of them except one (the starting node). The
value $C = C(\Delta)$ can be chosen appropriately, for any fixed
$\Delta$.} are
larger than $C(\Delta)$. For example, $\Delta= 3$ which
implies that there are at most $2 (\lfloor p /3 \rfloor+1)$ edges with
nonzero coefficients being smaller than $C(\Delta=3)$, and at least
$p- \lfloor
p/3 \rfloor-2$ edges with nonzero coefficients larger than
$C(\Delta=3)$. This implies that Condition
\ref{beta-mincondition}, given below, holds with a value $C(\Delta=3)$ of
order 1.

\section{Conditions and main result}\label{conditionsmainsection}

We write $\Sigma_0 =:( \sigma_{k,j } ) $, and we let
$\sigma_j^2:= \sigma_{j,j} $, $j=1,\ldots, p $.

\begin{condition}\label{sigmacondition} For some constant $\sigma
_0^2$, it
holds that
\[
\max_{1 \le j \le p } \sigma_j^2 \le
\sigma_0^2.
\]
\end{condition}

\begin{condition} \label{Lambda-mincondition} The smallest eigenvalue
$\Lambda_{\mathrm{min}}^2 $ of $\Sigma_0$ is nonzero; see also~(\ref{omega-pos}).
\end{condition}

\begin{condition}\label{hatedgescondition} For a given constant
$\alpha>0$,
it holds that for any $B= ( \beta_1,\ldots, \beta_p) \in{\cal B}$,
where ${\cal B}$ is as in (\ref{ell0MLE}), that
$s_{\beta_j} \le\alpha n / \log p $ for all $j=1,\ldots, p $, where
for a vector $\beta\in\R^p$ we denote the cardinality of its support
set by $s_{\beta}:= \# \{ \beta_k
\not= 0 \}$.
\end{condition}

\begin{condition} \label{edgescondition}
For some constant $\tilde\alpha$
and any permutation $\pi$, and all~$j$,
\[
\tilde s_j ( \pi) \le\tilde\alpha n / \log p,
\]
where $\tilde s_j ( \pi)=s_{\tilde\beta_j^0} $ is the number of incoming
edges of the DAG $(\tilde B_0 (\pi), \tilde\Omega_0 (\pi))$ at node
$j$; see also (\ref{cond-edges}).
\end{condition}

\begin{condition} \label{beta-mincondition}
There exist constants $0 \le\eta_1< 1 $ and $0 < \eta_0^2 < 1- \eta_1
$, such that for any permutation $\pi$,
the DAG $ ( \tilde B_0 (\pi), \tilde\Omega_0 ( \pi) ) $ [which has
$\tilde s( \pi) $ edges] has at least $(1- \eta_1) \tilde s (\pi)$
edges $(k,j)$ with $| \tilde\beta_{k,j}^0 (\pi) | > \sqrt{ \log p / n }
( \sqrt{p/ s_0 } \vee1 ) / \eta_0$.
\end{condition}

Following \citet{BvdG2011}, we refer to Condition~\ref{beta-mincondition}
as the ``beta-min'' condition. It is a ``kind of replacement'' of the
strong faithfulness condition in (\ref{strong-faith}) that is required for
consistency of the PC algorithm and variants thereof; see Section
\ref{subsecDAGSEM}. A detailed discussion about the assumptions is given
in Section~\ref{discussionsection}.

In the current section, we will present an asymptotic formulation for clarity.
We will provide a nonasymptotic result in Section~\ref{proofssection}.
Our results depend on $\Sigma_0$ via the constants
$\sigma_0$, $\Lambda_{\mathrm{min}}$, $\eta_1$ on the further constants
$\gamma_0:= (\alpha, \tilde\alpha, \eta_0)$
used in the conditions and on
$\alpha_0$ where
$1- \alpha_0$ is the confidence level of the statement.

We assume that we can take
$\gamma_0 $ sufficiently small.
Moreover, we state the results with $\alpha_0:= (4/p)\wedge0.05$ to avoid
digressions concerning the confidence level.
Explicit expressions can be found in Theorem~\ref{exact2theorem},
where we simplified the situation by assuming $n$ is sufficiently large and
$\log p / n $ is sufficiently small.

With the notation $z= {\mathcal O}(1)$
we mean that $z$ can be bounded by a constant depending only
on $\sigma_0 $ and $\Lambda_{\mathrm{min}}$.
Moreover,\vspace*{1pt} with
$z\asymp1$
we mean $z = {\mathcal O} (1)$ and $1/z= {\mathcal O} (1) $. Furthermore,
the Frobenius norm is defined as $\|B\|_F = (\sum_{j,k=1}^p
|\beta_{k,j}|^2)^{1/2}$ for a $p \times p$ matrix $B$ with elements
$\beta_{k,j}$.

\begin{theorem} \label{maintheorem} Assume Conditions \ref
{sigmacondition},~\ref{Lambda-mincondition},
\ref{hatedgescondition},~\ref{edgescondition} and \ref
{beta-mincondition}, with
$\gamma_0:=( \alpha, \tilde\alpha, \eta_0)$ sufficiently
small, but allowing
$1/ \| \gamma_0 \|_1 = {\mathcal O} (1) $.
Let $1-\alpha_0$ be the confidence level, with $\alpha_0:= (4
/p)\wedge0.05$.
Then for a choice
\[
\lambda^2 \asymp{\log p \over n} \biggl(
{ p \over s_0 } \vee1 \biggr),
\]
it holds that with probability at least $1- \alpha_0$,
\[
\| \hat B - \tilde B_0 \|_F^2 + \| \hat
\Omega- \tilde\Omega_0 \|_F^2 = {\mathcal O}
\bigl( \lambda^2 s_0\bigr),
\]
where $(\tilde{B}_0,\tilde{\Omega}_0)$ are defined in Section
\ref{subsecpermutations}, and
\[
\hat s \asymp\tilde s \asymp s_0.
\]
\end{theorem}

The proof is given in Section~\ref{proofssection}.
Theorem~\ref{exact2theorem} gives some explicit bounds.

\begin{remark}
Note that if the true permutation $\pi_0$ defined by
$\tilde B_0 (\pi_0) = B_0$ is known, then
from (multiple) regression theory the optimal rate of convergence
in Frobenius norm will be of order $(p+s_0) \log p / n $ (with $p \log
p / n $ being a lower bound
due to estimating the $p$ residual variances). Hence, Theorem \ref
{maintheorem} says that this
rate can also be achieved when not knowing $\pi_0$.
As in a multiple regression setup, a natural normalization of the
Frobenius norm is to divide it by $p$.
With this normalized norm, the estimator is consistent when the average
number of incoming edges $s_0 / p $ is of small order $n / \log p
$.\looseness=-1
\end{remark}

\begin{remark}\label{subsecexact}
If the beta-min condition (Condition~\ref{beta-mincondition}) holds with
$\eta_1=0$ and with very small values for $\eta_0:= \eta_0 (\pi)$,\vspace*{1pt} namely
of order $1/ \tilde s(\pi)$, then one obtains the screening\vspace*{2pt} property: all
edges in $( \tilde B_0, \tilde\Omega_0)$ are then with large probability
also present in $( \hat B, \hat\Omega)$.

Moreover, by taking $\lambda^2:= \lambda^2 (s_0)$ very
large (of order $s_0 \log p / n $), one can obtain with large
probability that\vadjust{\goodbreak}
$\hat s \le s_0$. In other words,
by imposing a strong beta-min condition, which is severe if $s_0$ is
large, one recovers with high probability the edges of the minimal-edge I-MAP
exactly. However, in Theorem~\ref{maintheorem}, we do not use
such values for $\eta_0$ and $\lambda$, but instead $\lambda^2 \asymp
\log p
/ n $ [when $p = {\mathcal O} (s_0)$] and
$\eta_0 \asymp1 $. Thus, we generally do not recover the true edges. This
is the price for dealing with
a large $p$ situation and an $s_0$ possibly growing in $n$. Such problems
do not show up in asymptotics with $p$ fixed.
\end{remark}

\section{A discussion of the conditions} \label{discussionsection}

\subsection{Bounds for the noise variances} \label{omegasection}
For all $\pi$ and $j$ and for any $\beta_j$ with
$\beta_{j,j} = 0$, we have $\| X_j - X \beta_j \|= \| X \beta_{j}^{-}
\| $ where
$\beta_{k,j}^{-} = -\beta_{k,j} $ for $k \not= j$ and \mbox{$\beta_{j,j}^{-}
= 1 $}.
It follows that for any $\pi$ and $j$,
\[
\bigl| \tilde\omega_j^0 (\pi) \bigr|^2 = \bigl\|
X_j - X \tilde\beta_j^0 (\pi)
\bigr\|^2 \ge\Lambda_{\mathrm{min}}^2
\]
with\vspace*{1pt} $\Lambda_{\mathrm{min}}^2$ the smallest eigenvalue of $\Sigma_0$.
Moreover, clearly
$|\tilde\omega_j^0 (\pi) |^2 \le\sigma_j^2 $. Hence, Conditions
\ref{sigmacondition} and~\ref{Lambda-mincondition} imply that for
all $\pi$ and $j$
\[
0< \Lambda_{\mathrm{min}}^2 \le\bigl|\tilde\omega_j^0
(\pi) \bigr|^2 \le\sigma_0^2.
\]

Furthermore, $\Lambda_{\mathrm{min}}^2 > 0$ is implied by
%
\begin{equation}
\label{omega-pos} \min_j \bigl|\omega_j^0\bigr|^2
> 0,
\end{equation}
since $\det(\Sigma_0) = \det(\Omega_0) = \prod_{j=1}^p \omega_j^0$. Thus,
Condition~\ref{Lambda-mincondition} is equivalent to (\ref{omega-pos}).

\subsection{Overfitting} \label{edgessection}
Condition~\ref{hatedgescondition} will ensure that the penalized
minus log-likelihood cannot become minus infinity.
If $n$ or more edges are allowed at a node, say at node $j$, the estimator
will overfit the data at this node, giving a residual variance $\hat
\omega_j^2= 0 $. The penalized minus log-likelihood is proportional to
$\sum_{j=1}^p \log\hat\omega_j^2 + \lambda^2 \hat s$ which will be
$- \infty$ if one allows that $\hat\omega_j$
vanishes. Note that the penalty as such does not prevent this type of
overfitting. Therefore, we need
a restriction on the class of possible DAGs, and Condition \ref
{hatedgescondition} serves this purpose.
We will show in Lemma~\ref{omegalemma} that Conditions \ref
{sigmacondition},~\ref{Lambda-mincondition} and \ref
{hatedgescondition} imply
that for an appropriate constant $K_0 >0$, it holds for all $j$ that
$\hat\omega_j \ge1/ K_0 $ with large probability.

\subsection{The beta-min condition} \label{beta-minsection}
One may circumvent the beta-min condition if one allows for edges with
weights below some noise level $\lambda_*$
to be set to zero. Here, $\lambda_*:= \sqrt{\log p / n} / \eta_0^*$
for some suitable $\eta_0^* > 0$.
Instead of trying to estimate the true DAG $(B_0, \Omega_0)$, one now
aims at estimating its best sparse approximation
$(B_0^*, \Omega_0^*) $, which is defined as follows.
Let for any DAG $(B, \Omega)$, and for $\Theta= \Theta( B, \Omega
)$, the weights $B_{\Theta} ( \pi)$ be obtained
by doing the Gram--Schmidt orthogonalization for $\| \cdot\|_{\Sigma
}$, where $\Sigma= \Theta^{-1}$
and $\| X \beta\|_{\Sigma}^2:= \beta^T \Sigma\beta$, $\beta\in\R
^p $.
Thus $B_{\Theta} (\pi) $ is lower-diagonal after the\vadjust{\goodbreak} permutation $\pi$
of its rows and columns, and for appropriate
$\Omega_{\Theta} (\pi)$, the DAG $( B_{\Theta} (\pi), \Omega_{\Theta}
(\pi)) $ satisfies
\[
\Theta= \Theta\bigl( B_{\Theta} (\pi), \Omega_{\Theta} (\pi)\bigr).
\]
Let $s_{\Theta} (\pi):= s_{B_{\Theta} ( \pi)} $ be the number of edges
of $B_{\Theta} (\pi) $.
Connecting this with our previous notation, we note that
\[
B_{\Theta_0} (\pi) = \tilde B_0 (\pi),\qquad \Omega_{\Theta_0} (
\pi) = \tilde\Omega_0 (\pi),\qquad s_{\Theta_0 } (\pi) = \tilde s ( \pi).
\]
Let now for some constant $\eta_0^*>0$,
\[
s_{\Theta}^* (\pi):= \# \bigl\{ (k,j)\dvtx  \bigl| \beta_{\Theta,k,j} (\pi) \bigr| >
\sqrt{\log p / n} \bigl( \sqrt{ p / s_{\Theta} (\pi) } \vee1 \bigr) /
\eta_0^* \bigr\}.
\]
We then take
\[
\Theta_0^*:= \arg\min\bigl\{ l(\Theta)\dvtx  \Theta= ( B, \Omega) \mbox{ a
DAG}, s_{\Theta}^* (\pi) \ge\bigl(1- \eta_1^*\bigr)
s_{\Theta} (\pi), \forall\pi\bigr\},
\]
where $0 \le\eta_1^* < 1$ and $l ( \Theta) = \operatorname{trace} (\Theta\Sigma
_0 ) - \log \operatorname{det} (\Theta)
= \EE l_n (\Theta)$ is the theoretical counterpart of the minus log-likelihood.
[Note that $\Theta_0:= \Sigma_0^{-1} $ is the overall minimizer of $l(
\Theta)$.]
We let $(B_0^*, \Omega_0^*)$ be a solution of
\[
\Theta_0^* = \Theta_0^* \bigl( B_0^*,
\Omega_0^* \bigr)
\]
with the minimum number of edges.
With constants $\eta_0^*$ and $\eta_1^*$ sufficiently small,
one may replace $\Theta_0=\Theta_0 ( B_0, \Omega_0) $ by $\Theta_0^*=
\Theta_0^* (B_0^*, \Omega_0^*)$ in our analysis.
In this way, one can avoid the beta-min condition, provided that the
bias term
that will now appear in the bounds is small enough.

\subsubsection{The beta-min condition and the number of edges}\label{subsecbetaminedges}
We further note that
Conditions~\ref{sigmacondition},~\ref{Lambda-mincondition} and
\ref{beta-mincondition} imply Condition~\ref{edgescondition} with
%
\begin{equation}
\label{cond-edges} \tilde{\alpha} = { \sigma_0^2 \eta_0^2 \over\Lambda
_{\mathrm{min}}^2 (1- \eta_1) }.
\end{equation}
This is because for all $j$,
\[
{ (1-\eta_1) \tilde s_j (\pi) \over\eta_0^2 } { \log p \over n} \le\bigl\|
\tilde
\beta_j^0 (\pi) \bigr\|_2^2 \le\bigl\| X
\tilde\beta_j^0 (\pi) \bigr\|^2 /
\Lambda_{\mathrm{min}}^2 \le\sigma_0^2 /
\Lambda_{\mathrm{min}}^2.
\]

\subsubsection{The strong beta-min condition in comparison to strong faithfulness}\label{subsecbetamin-strongfaith}

The beta-min condition as stated in Condition~\ref{beta-mincondition} is
of a rather weak form. In order to make a (vague)
relation to strong faithfulness, which focuses on exact edge recovery,
we consider the stronger version as discussed in Remark~\ref{subsecexact}.
As written in this remark, recovery of a minimal-edge I-MAP is
guaranteed with a
value for the lower bound on the weights of the nonzero edges of the order
$s_0 \sqrt{\log(p)/n} = p \sqrt{\log(p)/n}$ assuming $s_0 \asymp
p$: such a value in the beta-min condition is reasonable in the regime
$p =
o(\sqrt{n/\log(n)})$.

Although this seems rather restrictive at first sight, the
PC-algorithm necessarily requires restricted strong faithfulness
[\citet{uhleretal12} cf.] for consistent estimation of the Markov
equivalence class. Such a restricted strong faithfulness assumption has been
analyzed when assuming i.i.d. sampling of the nonzero edge weights. It
holds when assuming an upper bound on the growth of the dimension. The
best dimensionality range is achieved for bounded-degree
trees which restricts $p = o(\sqrt{n/\log(n)})$ to the same order of
magnitude as above while for other graphs the constraint on $p$ can be much
stronger, for example, $p = o(\log(n))$ for certain bipartite graphs
[\citet{uhleretal12}, Section 5.1].

The beta-min Condition~\ref{beta-mincondition} is not directly
comparable to the (restricted) strong faithfulness assumption. Therefore,
we cannot make a direct comparison between our penalized maximum
likelihood estimator and the PC-algorithm. The AR(1) model in Section
\ref{subsecAR1} is an example where the beta-min
condition holds with a value of order 1. We can extend the analysis to an
AR($k$) model with fixed $k \ge
2$: using analogous reasoning as for the AR(1) in Section
\ref{subsecAR1}, the
beta-min condition holds for a value of order 1. The theory from
\citet{uhleretal12} regarding strong faithfulness cannot be used for
AR(k) models since the corresponding edge
weights involve $k$ values which are the same throughout the whole
graph, that is, no i.i.d. sampling.

We will discuss in Section~\ref{equalvariancessection} the case where
the error variances are the same, that is, $\Omega_0 = \omega^0 I$. We then
only need a beta-min condition for the true underlying DAG instead of all
permutations (see the discussion after Theorem
\ref{equalvariancestheorem}). Thus, for the scenario $p =
o(\sqrt{n/\log(n)})$ in the
equal variance case, the beta-min condition is
very reasonable for any DAG. This is in sharp contrast to the constraint
arising from restricted strong faithfulness: if the underlying DAG is say
a certain bipartite graph, the corresponding dimensionality for
consistent edge recovery is very severe, see above.

Finally, we note that when focusing on bounding false positive
selections as in Theorem~\ref{maintheorem}, the $\ell_0$-penalized
MLE is
justified for the $p \gg n$ setting.

\subsection{\texorpdfstring{The high sparsity regime where $s_0\ll p$}{The high sparsity regime where s0 << p}}
The reason why we see a term $p / s_0 \vee1$ appearing in the tuning
parameter $\lambda^2$ (see Theorem~\ref{maintheorem})
and in
the beta-min condition (Condition~\ref{beta-mincondition}) is due to
the estimation
of the $p$ unknown variances, which gives a term of order $p \log p /
n $ in our bounds
for the squared Frobenius norm.
If $s_0 \ll p $, the true DAG has many disconnected components, and in
fact it then has
many isolated points. The variables in one component
are uncorrelated with those in another component. We see this in the zeroes
in the matrix $\Sigma_0$. The connected components and isolated points
are easily detected by $\Sigma_n$, assuming
that nonzero correlations are at least $ \sqrt{\log p / n }/ \eta_c$
in absolute value for an appropriate (sufficiently small)
constant $\eta_c$. Then we can do the analysis connected component by
connected component.
To summarize, the situation $p= {\mathcal O} (s_0)$ appears to be the most
interesting.
Alternatively, when the noise variances $\Omega_0$ are known up to a scalar
(e.g., if it is known that all noise variances are equal), we need not estimate
these variances anymore, and the term of order $p \log p / n$ does not
appear in the bounds, provided
an identifiability condition on the noise variances holds and $p $ is
sufficiently smaller than
$n/ \log n $.
This will be shown in the next section.

\section{The case of equal variances} \label{equalvariancessection}

Suppose that the noise variances $\{ | \omega_j^0 |^2 \}_{j=1}^p
$ are known up to a multiplicative
scalar. To simplify the exposition, let us assume that
\[
\omega_1^0 = \cdots= \omega_p^0
= 1.
\]

The $\ell_0$-penalized maximum likelihood estimator now becomes
%
\begin{equation}
\label{ell0eqvarMLE} \hat B := \arg\min \bigl\{\operatorname{trace} \bigl( ( I-B)
(I-B)^T \Sigma_n \bigr) + \lambda^2
s_B\dvtx
( B, I) \mbox{ a DAG}, B \in{\cal B} \bigr\},\hspace*{-35pt}
\end{equation}
where $B$ is as in (\ref{ell0MLE}).

The main Theorem~\ref{maintheorem} as well as the remarks in Section
\ref{subsecexact} apply to the estimator (\ref{ell0eqvarMLE}) as
well, assuming exactly the same Conditions~\ref{sigmacondition}--\ref
{beta-mincondition}.

For the case where $p = {\mathcal O}(n/\log(n))$ is sufficiently small, we
obtain consistent estimation
of the true underlying DAG, and we gain in comparison to the main Theorem
\ref{maintheorem} by excluding the additional factor $(p/s_0 \vee1)$. We
make the following assumptions.
%
\begin{condition}\label{omega-mincondition}
$\!\!\!$There exists a constant $\eta_{\omega} >0$ such that for all
\mbox{$\tilde\Omega_0 (\pi) \not= I$},
\[
{1 \over p} \sum_{j=1}^p
\bigl(\bigl|\tilde{\omega}_j^0 (\pi) \bigr|^2 - 1
\bigr)^2 > 1/ \eta_{\omega}.
\]
\end{condition}

\begin{condition}\label{psmallcondition}
There exists a constant $\alpha_*$ such that
\[
p \le\alpha_* n / \log n.
\]
\end{condition}

We call Condition~\ref{omega-mincondition} the ``omega-min''
condition. It
leads to
identification of the DAG with equal variances. Condition \ref
{psmallcondition} ensures that the
rate of convergence is fast enough to ensure that eventually we choose
the right
permutation. Note that it implies Conditions~\ref{hatedgescondition} and
\ref{edgescondition}, with $\alpha= \tilde\alpha= \alpha_*$.

Let $\pi_0$ be defined by $\tilde B_0 ( \pi_0) = B_0$. Since $B_0$ is
identifiable from the observational distribution ${\cal N}(0,\Sigma_0)$
[\citet{Peters2012}] (see also Section~\ref{subsecmainresults}), $\pi_0$
corresponds to the unique true ordering of the variables.

\begin{theorem}\label{equalvariancestheorem} Assume Conditions \ref
{sigmacondition} and
\ref{Lambda-mincondition} and Conditions
\ref{omega-mincondition}\break and~\ref{psmallcondition}. Let $\alpha_0:=
(4/p )\wedge0.05$.
Then for $\gamma_*:= (\alpha_*, \eta_{\omega})$ suitably small, but
allowing
$1/ \| \gamma_* \|_1 = {\mathcal O} (1)$, and for
$\lambda^2 \asymp\log p / n $, it holds
with probability at least
$1- \alpha_0$, that
$\hat\pi= \pi_0$,
and
\[
\| \hat B - B_0 \|_F^2 +
\lambda^2 \hat s = {\mathcal O} \bigl( \lambda^2
s_0 \bigr).
\]
\end{theorem}

The proof is given in Section~\ref{equalvariancesproofsection}.\vadjust{\goodbreak}

Thus, we find $\hat s = {\mathcal O} (s_0)$, but we do not show $\hat s
\asymp s_0$. To establish the latter, one again needs a beta-min condition,
but this time only on the DAG $( B_0, I)$, and not on any of the other
representations $(\tilde B_0 (\pi), \tilde\Omega_0 (\pi) )$ with $\pi
\not= \pi_0$. This is a much simplified and weaker assumption than in Condition
\ref{beta-mincondition}. Furthermore, since $\hat{\pi} = \pi_0$ with large
probability, a refit of the model using a (de-coupled) penalized node-wise
regression with parents according to $\hat{\pi}$ will with large
probability recover the edges under the standard conditions for such a
method [e.g., a node-wise Lasso will with large probability recover the edges
under the condition that for all $j$, $|\beta^0_{k,j}| >
\sqrt{\tilde{s}_j(\pi_0) \log(p)/n}/\eta^0$ for some sufficiently small
$\eta^0 > 0$].

\subsection{The non-Gaussian case}
To avoid technical digressions in our proofs, we assume a
Gaussian distribution for the observations where zero correlations mean
independence. We
use in Lemma~\ref{Gaussianlemma} that if for some $\tilde\epsilon_j$,
$\EE\tilde\epsilon_j=0$, then also the conditional expectation of
$\tilde
\epsilon_j$ given variables $X_k$ that are
uncorrelated with $\tilde\epsilon_j$ is zero. In the non-Gaussian case,
this is no longer true.
However, one can still derive similar results, along a line of proof that
does not use conditioning but instead concentration inequalities for
averages of products of random
variables (empirical covariances). This means that our results go through
for observations which are sub-Gaussian. The proofs then rely on concentration
inequalities of Bernstein-type. The main adjustments of our proofs are then
as follows. We assume that the rows of $X$ form an i.i.d.
sequence of sub-Gaussian vectors as defined in \citet{loh2012} and
replace Theorem~\ref{normtheorem} by their Lemma 15.
In Lemma~\ref{Gaussianlemma} we assume $\epsilon$ and $Z$ are
sub-Gaussian and uncorrelated, and replace the empirical
squared norm $\| X \beta\|_n^2:= [ \sum_{i=1}^n ( X \beta)_i^2 /n ]$
by the theoretical squared norm
$\| X \beta\|^2 $. We can then apply similar arguments as used in
Lemma 15 of
\citet{loh2012}.
In Lemma~\ref{startlemma}, we no longer use the empirical squared norm
but instead the theoretical one. Theorem~\ref{T2theorem} needs virtually
no adjustments.

\section{Conclusions}
We establish the first results of the $\ell_0$-penalized MLE for estimation
of the minimal-edge I-MAP (the smallest DAG which can generate the
data-generating
distribution) in the high-dimensional sparse setting. Thereby, we avoid the
faithfulness condition, and the strong
faithfulness assumption (\ref{strong-faith}) or its restricted
version; cf. \citet{uhleretal12}; the latter is necessary for consistency
of the PC-algorithm [\citet{sgs00}]. The (restricted) strong faithfulness
condition
is typically very strong [\citet{uhleretal12}] and hence, our results
contribute in relaxing such very restrictive assumptions.

Our main assumption is
the beta-min Condition~\ref{beta-mincondition} (which implies the sparsity
Condition~\ref{edgescondition}; see Section~\ref{subsecbetaminedges}):
as an example, the AR(1)-model in Section~\ref{subsecAR1} fulfills
it, even if $p \gg n$. The noise level is of the
order $\sqrt{\log(p)/n} (p/s_0 \vee1)$: the additional factor $(p/s_0
\vee
1)$ occurs due to estimation of $p$ variances in $\Omega_0$. However, the
interesting scenario is for the case where $s_0 \ge \operatorname{const.} p$ since $s_0
\ll p$ corresponds to a DAG where most nodes are isolated having no
edges to
other nodes; thus, for $s_0 \ge \operatorname{const.}p$, we obtain the usual noise
level of the order $\sqrt{\log(p)/n}$, as in high-dimensional regression
problems.

For the equal variance case with $p = {\mathcal O}(n / \log(n))$
sufficiently small, our result in Theorem
\ref{equalvariancestheorem} (and its comment below) is most clear in that
we essentially only require the beta-min Condition
\ref{beta-mincondition} for the true DAG $B_0$, that is, a substantially
relaxed assumption, and the identifiability
Condition~\ref{omega-mincondition} for the error variances: we can then
recover the true underlying unique DAG $B_0$.
Thus, we have identified an important class of models
where estimation of the order of variables and the true underlying DAG is
possible without requiring the badly limiting (restricted) strong faithfulness
condition (\ref{strong-faith}).

\section{Proofs}\label{proofssection}

\subsection{A brief outline of the proofs}\label{proofsoutlinesection}

{We first consider the proof of Theorem~\ref{maintheorem} which treats
the case of unknown variances $\{ (\omega_j^0 )^2 \} $.
In Lemma~\ref{startlemma} of Section~\ref{subsetprobabilitysection}, we present a bound for $\sum_{j=1}^p [ ( \tilde
\omega_j^0)^2 /
\hat\omega_j^2 -1 ]^2 $ and for
$\sum_{j=1}^p \| X ( \hat\beta_j - \tilde\beta_j^0 )\|_n^2 $ using
the empirical norm $\| v \|_n:=$ $[ \sum_{i=1}^n v_i^2/n ]^{1/2} $, $v
\in\R^n $. The result follows from a straightforward manipulation of
likelihoods,
but it is assumed there that
one is on the part of the probability space where the random
components behave well.
The study of these random components
is deferred to Sections~\ref{T1section},~\ref{T2section} and~\ref
{T3T0section}. First, the bound of Lemma~\ref{startlemma} is refined
because it involves the
number of edges $\tilde s$ of the DAG formed by using the random
permutation $\hat\pi$ that the penalized maximum
likelihood estimator chooses. Section~\ref{exploitsection} presents
the tools to deal with this by exploiting the beta-min condition. The
idea here
is that if the Frobenius norm between $\hat B$ and $\tilde B$ is small,
the number of edges of $\tilde B$ cannot be much larger than
those of $\hat B$.

A substantial part of the proof of Theorem~\ref{maintheorem} goes into
showing that with large probability we are on a set of the form
$\bigcap_{k=0}^3{ \cal T}_k $ where the random components behave well.
Let us first discuss ${\cal T}_1$. Here, a uniform inequality
holds for the empirical correlation between the projections and error
terms in a Gram--Schmidt orthogonalization.
For a fixed permutation $\pi$ it is rather standard to control these
empirical correlations. The new element is that we have to control them
uniformly
over all permutations $\pi$ in order to show that ${\cal T}_1$ has
large probability. We do this in Section~\ref{T1section}, where the
arguments
used are explained just before Theorem~\ref{T1theorem}.
In the set ${\cal T}_2 $ all empirical variances of the error terms in
a Gram--Schmidt orthogonalization are
close to their expectations. We show in Section~\ref{T2section}
that uniformly over all $\pi$ this is true with large probability.
The set ${\cal T}_3$ gives bounds for $\| \beta\|_2$ in terms of $\|X
\beta\|_n$ and the number of nonzero
coefficients in $\beta$. We show in Theorem~\ref{normtheorem} that
${\cal T}_3$ has large probability.
This makes it possible to move from empirical norms to Frobenius norms
and moreover
shows that with large probability the $\{ \hat\omega_j^2 \}$ are bounded
away from zero. The latter event is defined as the set ${\cal T}_0$.

For the proof\vspace*{-2pt} of Theorem~\ref{equalvariancestheorem} where the
variances $(\omega_j^0 )^2 $ are all known to be equal to one, we
use the same structure. We assume that we are on the set $\bigcap_{k=1}^3
{\cal T}_k$, and use straightforward manipulations of likelihoods.

\subsection{Bounds on a subset of the probability space} \label{subsetprobabilitysection}
We present some explicit bounds assuming we are on a set of the form
$\bigcap_{k=0}^3 {\cal T}_k $, where the sets ${\cal T}_k$ are defined below.
Then we show in Sections~\ref{T1section},~\ref{T2section} and
\ref{T3T0section} that each ${\cal T}_k $, $k=0,\ldots, 3$, has large
probability
for an appropriate choice of the constants and of the parameters
$\lambda_1 $, $\lambda_2$ and $\lambda_3$ involved in the definition of these
sets. In fact, we will show that one can take
$\lambda_1 \asymp\lambda_2 \asymp\lambda_3 \asymp\sqrt{\log p / n} $.

Let for some constant $K_0>0$,
\[
{\cal T}_0:= \bigl\{
\hat\omega_j^2
\wedge 1/
\hat\omega_j^2
\ge1/K_0^2, \forall j \bigr\}.
\]
Let us write $X_k \perp\tilde\epsilon_j$ if $X_k$ and $\tilde\epsilon_j$
are independent. For all $\pi$
and $j$, define $\tilde\epsilon_j (\pi) = X_j - X \tilde\beta_{j}^0
(\pi) $, and
$\tilde{\cal B}_j (\pi):= \{ \beta_j\dvtx  X_k \perp\tilde\epsilon_j
(\pi),
\forall\beta_{k,j} \not= 0 \} $. Moreover, let $\tilde{\cal B} (
\pi):=
\{ B= ( \beta_1,\ldots, \beta_p) \in{\cal B}\dvtx  \beta_j \in
\tilde{\cal B}_j (\pi)\ \forall j \} $.
For
some $\delta_1 >0$ and some $\lambda_1 > 0 $, write
\begin{eqnarray*}
{\cal T}_1&:=& \Biggl\{2 \sum_{j=1}^p
\bigl| \tilde\epsilon_j^T(\pi) X\bigl( \beta_j -
\tilde\beta_j^0 (\pi) \bigr) \bigr| /n \\
&&\hspace*{6pt}\le
\delta_1 \sum_{j=1}^p \bigl\| X
\bigl( \beta_j - \tilde\beta_j^0 (\pi)
\bigr) \bigr\|_n^2
\\
&&\hspace*{14.6pt}{}+ \lambda_1^2 \bigl(s + \tilde s(\pi) \bigr)/
\delta_1, \forall B= (\beta_1,\ldots,
\beta_p) \in\tilde{\cal B} (\pi)\ \forall\pi\Biggr\}.
\end{eqnarray*}
We let for some $\lambda_2>0$,
\[
{\cal T}_2:= \Biggl\{ \sum_{j=1}^p
\biggl( { \| \tilde\epsilon_j
(\pi) \|_n^2 - | \tilde\omega_j^0 (\pi) |^2 \over| \tilde\omega_j^0
(\pi) |^2 } \biggr)^2 \le4 \lambda_2^2
\bigl(p+\tilde s (\pi) \bigr), \forall\pi\Biggr\},
\]
where we recall the notation $\| v \|_n^2:= v^T v / n, v \in\R^n $.
Finally, for some \mbox{$\delta_3>0$} and some $\lambda_3 > 0 $, let
${\cal T}_3$ be the set
\[
{\cal T}_3:= \bigl\{ \| X \beta\|_n \ge[
\delta_3 - \lambda_3\sqrt{ s_{\beta} } ] \| \beta
\|_2, \forall\beta\bigr\}.
\]
Recall that $s_{\beta}:= \# \{ \beta_{k} \not= 0 \} $.

\begin{lemma} \label{startlemma} Define $( \tilde B_0, \tilde\Omega
_0 ):=
( B_0 ( \hat\pi), \Omega_0 ( \hat\pi) ) $ and $\tilde s:=
s_{\tilde B_0 } $.
Assume that
Condition~\ref{sigmacondition} holds.
Suppose we are on $\bigcap_{k=0}^2 {\cal T}_k $ with
$0 < \delta_1 < 1/ K_0^2 $ and $0< \delta_2 < 1/ ( 2 K_0^4 \sigma_0^4 ) $.
Take the tuning parameter $\lambda^2 > \lambda_1^2 / \delta_1+ \lambda
_2^2 /\delta_2 $.
Then
\begin{eqnarray*}
&&
\biggl({1 \over K_0^2 } - \delta_1 \biggr) \sum
_{j=1}^p { \bigl\| X \bigl( \hat\beta_j -
\tilde\beta_j^0 \bigr) \bigr\|_n^2 } +
\biggl( { 1 \over2 K_0^4 \sigma_0^4 }- \delta_2 \biggr) \sum
_{j=1}^p \biggl( { \hat\omega_j^2 - | \tilde\omega_j^0| ^2 \over|
\hat\omega_j |^2 }
\biggr)^2
\\
&&\quad{}+ \biggl( \lambda^2 - { \lambda_1^2 \over\delta_1 }-
{ \lambda_2^2
\over\delta_2 } \biggr) \hat s \\
&&\qquad\le\lambda^2 s_0 +
{\lambda_2^2 (p+ \tilde s ) \over\delta_2 } + {
\lambda_1^2 \tilde s \over\delta_1 }.
\end{eqnarray*}
\end{lemma}

\begin{pf}
Let $\tilde\epsilon:= \tilde\epsilon(\hat\pi)$.
We apply the basic inequality
\[
l_n ( \hat\Theta) + \lambda^2 \hat s \le l_n
( \Theta_0 ) + \lambda^2 s_0
\]
or equivalently
\[
p + \sum_{j=1}^p \log\hat
\omega_j^2 + \lambda^2 \hat s \le\sum
_{j=1}^p { \| \epsilon_j \|_n^2 \over
| \omega_j^0 |^2 } + \sum
_{j=1}^p \log\bigl| \omega_j^0
\bigr|^2 + \lambda^2 s_0,
\]
which gives, using $\log(\operatorname{det}( \Sigma_0)) = \sum_{j=1}^p \log|
\omega_j^0 |^2 = \sum_{j=1}^p \log| \tilde\omega_j^0 |^2$,
\[
\sum_{j=1}^p \log\biggl(
{ \hat\omega_j^2 \over| \tilde\omega_j^0
|^2 } \biggr) + \lambda^2 \hat s \le\sum
_{j=1}^p \biggl( { \| \epsilon_j \|_n^2 \over
| \omega_j^0 |^2 } -1 \biggr)
+ \lambda^2 s_0.
\]
Since $\hat\omega_j^2 \ge1/K_0^2 $ (since we are on ${\cal T}_0$) and
$| \tilde\omega_j^0 |^2 \le\sigma_0^2 $
(by Condition~\ref{sigmacondition}), we know
that ${| \tilde\omega_j^0 |^2 / \hat\omega_j^2 } \le K_0^2 \sigma_0^2 $.
But then, using $ \log(1+x ) \le x - x^2 /(2 (1+ c )^2) $, $ -1 < x
\le c $,
we get
\[
\log\biggl( { \hat\omega_j^2 \over| \tilde\omega_j^0 |^2 } \biggr) = -
\log\biggl(
{ | \tilde\omega_j^0| ^2 \over| \hat\omega_j |^2 } \biggr) \ge- \biggl(
{ | \tilde\omega_j^0| ^2 \over| \hat\omega
_j |^2 } -1 \biggr)
+ { 1 \over2 K_0^4 \sigma_0^4 } \biggl( { | \tilde\omega_j^0| ^2 \over
| \hat\omega_j |^2 } -1
\biggr)^2.
\]
We plug this back into the basic inequality to get
\[
\sum_{j=1}^p { \hat\omega_j^2 - | \tilde\omega_j^0| ^2 \over| \hat
\omega_j |^2 } +
{ 1 \over2 K_0^4 \sigma_0^4 } \biggl( { \hat\omega_j^2 - | \tilde
\omega_j^0| ^2 \over| \hat\omega_j |^2 } \biggr)^2 +
\lambda^2 \hat s \le\sum_{j=1}^p
\biggl( { \| \epsilon_j \|_n^2
\over
| \omega_j^0 |^2 } -1 \biggr) + \lambda^2
s_0.
\]
Rewrite this to
\begin{eqnarray*}
&&
\sum_{j=1}^p { \| X ( \hat\beta_j - \tilde\beta_j^0 ) \|_n^2 \over
\hat\omega_j^2 } +
{ 1 \over2 K_0^4 \sigma_0^4 } \sum_{j=1}^p
\biggl( { \hat\omega_j^2 - | \tilde\omega_j^0| ^2 \over| \hat\omega_j
|^2 } \biggr)^2 + \lambda^2
\hat s
\\
&&\qquad
\le2 \sum_{j=1}^p { \tilde\epsilon_j^T X ( \hat\beta_j - \tilde\beta
_j^0 ) /n \over
\hat\omega_j^2 }
+ \sum_{j=1}^p \biggl(
{ \| \epsilon_j \|_n^2 \over
| \omega_j^0 |^2 } -1 \biggr) \\
&&\qquad\quad{}- \sum_{j=1}^p
\biggl( { \| \tilde
\epsilon_j \|_n^2 - | \tilde\omega_j^0 |^2 \over
| \hat\omega_j |^2 } \biggr) + \lambda^2 s_0.
\end{eqnarray*}

We now apply
\begin{eqnarray*}
&&
\sum_{j=1}^p \biggl(
{ \| \tilde\epsilon_j \|_n^2 - | \tilde\omega
_j^0 |^2 \over
| \hat\omega_j |^2 } \biggr)
\\
&&\qquad
= \sum_{j=1}^p \biggl(
{ \| \tilde\epsilon_j \|_n^2 - | \tilde\omega
_j^0 |^2 \over
| \tilde\omega_j^0 |^2 } \biggr) + \sum_{j=1}^p
\biggl( { \| \tilde
\epsilon_j \|_n^2 - | \tilde\omega_j^0 |^2 \over| \tilde\omega_j^0 |^2}
\biggr) \biggl( { | \tilde\omega_j^0 |^2 - \hat\omega_j^2 \over\hat
\omega_j^2 }
\biggr).
\end{eqnarray*}
But, by the Cauchy--Schwarz inequality and using that we are on ${\cal T}_2$,
\begin{eqnarray*}
&&
\Biggl| \sum_{j=1}^p \biggl(
{ \| \tilde\epsilon_j \|_n^2 - | \tilde\omega
_j^0 |^2 \over| \tilde\omega_j^0 |^2} \biggr) \biggl( { | \tilde\omega
_j^0 |^2 - \hat\omega_j^2 \over\hat\omega_j^2 } \biggr) \Biggr|
\\
&&\qquad
\le\Biggl( \sum_{j=1}^p \biggl(
{ \| \tilde\epsilon_j \|_n^2 - |
\tilde\omega_j^0 |^2 \over| \tilde\omega_j^0 |^2} \biggr)^2 \Biggr
)^{1/2} \Biggl( \sum
_{j=1}^p \biggl( { | \tilde\omega_j^0 |^2 - \hat\omega_j^2 \over\hat
\omega_j^2 }
\biggr)^2 \Biggr)^{1/2}
\\
&&\qquad\le2 \sqrt{ (p+\tilde s ) \lambda_2^2 } \Biggl(
\sum_{j=1}^p \biggl( { | \tilde\omega_j^0 |^2 - \hat\omega_j^2 \over\hat
\omega_j^2 }
\biggr)^2 \Biggr)^{1/2} \\
&&\qquad\le{ (p+\tilde s ) \lambda_2^2 \over\delta_2 } +
\delta_2 \sum_{j=1}^p \biggl(
{ | \tilde\omega_j^0 |^2 - \hat\omega_j^2 \over\hat\omega_j^2 } \biggr)^2.
\end{eqnarray*}

Invoking $\operatorname{trace}( \Theta_0 \Sigma_n ) = \operatorname{trace} (\tilde\Theta
_0 \Sigma_n) $, that is,
\[
\sum_{j=1}^p { \| \epsilon_j \|_n^2 / \bigl| \omega_j^0 \bigr|^2 } =
\sum_{j=1}^p { \| \tilde\epsilon_j \|_n^2 / \bigl| \tilde\omega_j^0 \bigr|^2 },
\]
and using that we are on ${\cal T}_1$, we see that
\begin{eqnarray*}
&&
\biggl({1\over K_0^2 } - \delta_1 \biggr) \sum
_{j=1}^p { \bigl\| X \bigl( \hat\beta_j -
\tilde\beta_j^0 \bigr) \bigr\|_n^2 }\\
&&\quad{} +
\biggl( { 1 \over2 K_0^4 \sigma_0^4 }- \delta_2 \biggr) \sum
_{j=1}^p \biggl( { \hat\omega_j^2 - | \tilde\omega_j^0| ^2 \over|
\hat\omega_j |^2 }
\biggr)^2
+ \biggl( \lambda^2 - { \lambda_1^2 \over\delta_1 } \biggr) \hat s \\
&&\qquad\le
\lambda^2 s_0 + { \lambda_2^2(p+\tilde s ) \over\delta_2 } +
{
\lambda_1^2 \tilde s\over\delta_1 }.
\end{eqnarray*}
\upqed\end{pf}

\subsection{Exploiting the beta-min condition}\label{exploitsection}

\begin{lemma}\label{beta-minlemma}
Let $\tilde s = s_{\tilde B}$ be the number of edges of $\tilde B_0$
and $\hat s= s_{\tilde B_0 } $ be the number of
edges of $\hat B$. Suppose that for some
$\tilde\lambda$,
\[
\| \hat B - \tilde B_0 \|_F \le\tilde\lambda\sqrt{
\tilde s}
\]
and
that for some constant $0\le\eta_1< 1 $ and $0 < \eta_2^2 < 1- \eta_1$
\[
\# \bigl\{ \bigl| \tilde\beta_{j,k}^0 \bigr| \ge\tilde\lambda/
\eta_2 \bigr\} \ge(1- \eta_1 ) \tilde s.
\]
Then $\hat s \ge(1- \eta_1- \eta_2^2) \tilde s $.
\end{lemma}

\begin{pf}
Let
\[
{\cal N}:= \bigl\{(k,j)\dvtx  \bigl| \tilde\beta_{k,j}^0 \bigr| \ge\tilde
\lambda/ \eta_2 \bigr\},\qquad {\cal M}:= \bigl\{ (k,j)\dvtx  \bigl|\hat
\beta_{k,j} - \tilde\beta_{k,j}^0 \bigr| \ge\tilde
\lambda/\eta_2 \bigr\}.
\]
Then for $(k,j) \in{\cal N} \cap{\cal M}^c$ it holds that
\[
| \hat\beta_{k,j} | \ge\bigl| \tilde\beta_{k,j}^0 \bigr| - \bigl|
\hat\beta_{k,j} - \tilde\beta_{k,j}^0 \bigr| > 0,
\]
so that $\hat s \ge| {\cal N} \cap{\cal M}^c | $.
Since $\| \hat B - \tilde B_0 \|_F \le\tilde\lambda\sqrt{\tilde s}
$, we must have
\[
\sum_{(k,j) \in{\cal N} \cap{\cal M} } \bigl| \hat\beta_{k,j} - \tilde
\beta_{k,j}^0 \bigr|^2 \le\sum
_{(k,j) } \bigl| \hat\beta_{k,j} - \tilde\beta_{k,j}^0
\bigr|^2 = \| \hat B - \tilde B_0 \|_F^2
\le\tilde\lambda^2 \tilde s,
\]
whereas
\[
\sum_{(k,j) \in{\cal N} \cap{\cal M} } \bigl| \hat\beta_{k,j} - \tilde
\beta_{k,j}^0 \bigr|^2 \ge| {\cal N} \cap{\cal M} |
\tilde\lambda^2/\eta_2^2.
\]
Hence $ | {\cal N} \cap{\cal M} | \le\eta_2^2 \tilde s $.
This gives
\[
\bigl| {\cal N} \cap{\cal M}^c \bigr| = | {\cal N}| - | {\cal N} \cap{\cal M}
| \ge(1- \eta_1) \tilde s - \eta_2^2 \tilde
s = \bigl(1- \eta_1- \eta_2^2\bigr) \tilde
s.
\]
\upqed\end{pf}

\begin{lemma}\label{beta-min2lemma}
Suppose that for some $\delta_B>0$, $\delta_s >0$, $\lambda_0>0$ and
$\lambda$ one has
\[
\delta_B \|\hat B - \tilde B_0 \|_{F}^2
+ \lambda^2 \delta_s \hat s \le\lambda^2
s_0 + \lambda_0^2 \tilde s,
\]
where $\tilde s \ge s_0$. Let $\tilde\lambda^2 \delta_B \ge\lambda^2
+ \lambda_0^2 $ and
assume
that
\[
\# \bigl\{ \bigl| \tilde\beta_{j,k}^0 \bigr| \ge\tilde\lambda/
\eta_2 \bigr\} \ge(1- \eta_1) \tilde s.
\]
Then
\[
\delta_B \|\hat B - \tilde B_0 \|_{F}^2
+ \biggl( \lambda^2 \delta_s - { \lambda_0^2 \over1- \eta_1^2 - \eta_2^2}
\biggr) \hat s \le\lambda^2 s_0
\]
and $\hat s \ge(1- \eta_1 - \eta_2^2 ) s_0 $.
\end{lemma}

\begin{pf}
Since $\tilde s \ge s_0 $, we find that
\[
\delta_B \| \hat B - \tilde B_0 \|_F^2
\le\bigl(\lambda^2+ \lambda_0^2 \bigr)\tilde
s \le\delta_B \tilde\lambda^2 \tilde s.
\]
This gives by Lemma~\ref{beta-minlemma} that $\hat s \ge(1- \eta_1 -
\eta_2^2 ) \tilde s $.
But then
\[
\delta_B \|\hat B - \tilde B_0 \|_{F}^2
+ \biggl( \lambda^2 \delta_s - {\lambda_0^2 \over1- \eta_1 - \eta_2^2 }
\biggr) \hat s \le\lambda^2s_0.
\]
\upqed\end{pf}

\subsection{The sets ${\cal T}_k$, $k=0,1,2,3$} \label{T-section}

\subsubsection{The set ${\cal T}_1$}\label{T1section}

\begin{lemma} \label{Gaussianlemma}
Let $Z$ be a fixed $n \times m$ matrix and $\varepsilon_1,\ldots,
\varepsilon_n$
be independent ${\cal N} (0, \sigma_0^2)$-distributed random variables.
Then for all $t > 0 $
\[
\PP\Bigl( \sup_{\| Z \beta\|_n \le1 } \bigl| \varepsilon^T Z \beta\bigr| / n
\ge\sigma_0 ( \sqrt{ 2m / n } + \sqrt{2 t / n} ) \Bigr) \le\exp[-t].
\]
\end{lemma}

\begin{pf}
Assume without loss of generality that $Z^T Z/n =I$ and
define $V_k:= \varepsilon^T Z_k /(\sigma_0 \sqrt n) $. Then
$V_1,\ldots, V_p$ are independent and ${\cal N} (0, 1)$-distributed.
It follows that for all $N \in
\{ 2,3, \ldots\} $, that $\EE|V_k^2| ^{N} = {(2N)! / ( 2^N N!)} \le
N! $.
But then by Bernstein's inequality [see \citet{Bennet62}], for all $t>0$,
%
\begin{equation}
\label{Bernstein} \PP\Biggl( \sum_{k=1}^m
\bigl( V_k^2 - \EE V_k^2 \bigr)
\ge2 \sqrt{t m} + 2 t \Biggr) \le\exp[-t].
\end{equation}
Now use that $\sum_{j=1}^m \EE V_k^2 = m$. We get
\[
\PP\Biggl( \sum_{k=1}^m
V_k^2 \ge m + 2 \sqrt{t m} + 2t \Biggr) \le\exp[-t].
\]
But $m + 2 \sqrt{tm} + 2 t \le( \sqrt{2m} + \sqrt{2 t})^2 $.
Furthermore,
\[
\sup_{\| Z \beta\|_n \le1 } \bigl| \varepsilon^T Z \beta\bigr| / n =
{ \sigma_0 \over\sqrt n } \sqrt{ \sum_{k=1}^m
V_k^2 /n }.
\]
\upqed\end{pf}

We are\vspace*{1pt} dealing now with the problem of uniformly controlling over all
permutations $\pi$.
We consider the local structure at each node of a DAG $(\tilde B_0 (\pi
), \tilde\Omega_0 (\pi)) $ with
$\tilde B_0 (\pi) =: (\tilde\beta_{k,j}^0 ( \pi) )$. Let
$\tilde S_j ( \pi)$ be the set of incoming edges at node $j$. Given
$\tilde S_j (\pi)$, the vector $X \tilde\beta_j^0(\pi) $ is the
projection in $L_2 (P)$ of $X_j $ on the linear space spanned by $\{
X_k \}_{k \in\tilde S_j (\pi) } $. Moreover,
$\tilde\epsilon_j ( \pi)$ is the anti-projection $\tilde\epsilon_j
(\pi) = X_j - X \tilde\beta_{j}^0 (\pi) $. In other words (for $j$ fixed)
if the parents $\tilde S_j (\pi) $ at node $j$ are given, then the
coefficients $\tilde\beta_{k,j}^0(\pi) $ and noise term $\tilde
\epsilon_j (\pi) $ are given as well.
Also, the set of variables $X_k$ that are independent of $\tilde
\epsilon_j$ is then given.
Recall that
$\tilde{\cal B}_j (\pi):= \{ \beta_j\dvtx  X_k \perp\tilde\epsilon_j
(\pi),
\forall\beta_{k,j} \not= 0 \} $.
Thus,\vspace*{1pt} for each fixed $j$, if $\tilde S_j (\pi)$ is given then the
local situation
$( \tilde\epsilon_j (\pi), \tilde
\beta_j^0 (\pi), \tilde{\cal B}_j (\pi) )$ at node $j$ is given.

Let $\Pi_j (m)$ be the collection of all permutations giving DAGs
$(\tilde B_0 (\pi), \tilde\Omega_0 (\pi))$
with edges
$( \tilde S_1 (\pi),\ldots, \tilde S_p (\pi) )$ with $|\tilde S_j
(\pi) | = m$.
If for some $m \in\{ 0, 1,\ldots, p \}$, we know that $\pi\in\Pi
_j(m)$; that is,
we know that node $j$ has $m$ parents, so then there are at most ${ p
\choose m} $
possibilities for the local situation at node $j$.

\begin{theorem} \label{T1theorem} Assume Condition~\ref{sigmacondition}.
Then for all $t > 0$,
\begin{eqnarray*}
&&
\PP\Biggl(\max_{\pi} \sup_{B \in\tilde{\cal B} (\pi) } { 2 \sum
_{j=1}^p \bigl| \tilde\epsilon_j^T
(\pi) X \bigl( \beta_j - \tilde\beta_j^0
(\pi) \bigr) \bigr| / n - \delta_1 \sum_{j=1}^p
\bigl\| X \bigl( \beta_j - \tilde\beta_j^0 (\pi)
\bigr) \bigr\| _n^2 }
\\
&&\qquad\hspace*{85.5pt}
\ge{ 4 \sigma_0^2 ( s_B+ \tilde s (\pi) ) \over n \delta_1 } + { \sigma
_0^2 (t+ 2\log p )(s_B + \tilde s(\pi) ) \over n \delta_1 } \Biggr)
\\
&&\qquad\le\exp[ - t ].
\end{eqnarray*}
\end{theorem}

\begin{pf}
Let $A_j (\pi)$ be the event
\begin{eqnarray*}
A_j (\pi) &:=& \biggl\{ \exists\beta_j \in\tilde{\cal
B}_j (\pi)\dvtx  \sup_{\| X ( \beta_j - \tilde\beta_j^0 (\pi)) \|_n \le1 } \bigl|
\tilde
\epsilon_j^T (\pi) X \bigl( \beta_j -
\tilde\beta_j^0 (\pi)\bigr) \bigr| / n
\\
&&\hspace*{6.5pt}
\ge\sigma_0 \biggl( \sqrt{ 2 (s_{\beta_j} + \tilde s_j (\pi) ) \over n
} + \sqrt
{2 (t + \tilde s_j (\pi) \log p + 2 \log p) \over n } \biggr) \biggr\}.
\end{eqnarray*}
Then
by Lemma~\ref{Gaussianlemma}, for all $t>0$, $\pi$ and $j$
\[
\PP\bigl( A_j (\pi) \bigr) \le\exp\bigl[ -\bigl(t+ \tilde
s_j (\pi) \log p + 2\log p \bigr) \bigr].
\]
We now let $\pi$ vary over all permutations such that $|\tilde S_j (\pi
) | = m$.
We then get
\[
\PP\biggl( \bigcup_{\pi\in\Pi_j (m ) } A_j (\pi)
\biggr) \le\pmatrix{ p \cr m} \exp\bigl[ -(t+ m \log p + 2\log p ) \bigr]
\le\exp
\bigl[-(t + 2 \log p ) \bigr].
\]
Next, we let $\pi$ vary over all permutations. We get
\begin{eqnarray*}
\PP\biggl( \bigcup_{\pi} A_j (\pi)
\biggr) &\le&\sum_{m=1}^p \max
_{1 \le
m \le p } \PP\biggl( \bigcup_{\pi\in\Pi_j (m) }
A_j (\pi) \biggr)
\\
&\le& p \exp\bigl[-(t + 2 \log p )\bigr]\\
&\le&\exp\bigl[-(t+ \log p ) \bigr].
\end{eqnarray*}
Finally
\[
\PP\Biggl( \bigcup_{j=1}^p \bigcup
_{\pi} A_j (\pi) \Biggr) \le p \max
_{j} \PP\biggl( \bigcup_{\pi
}
A_j (\pi)\biggr) \le p \exp\bigl[-(t+\log p ) \bigr] \le\exp[-t].
\]

Now, we use that for all $\delta_1 > 0$,
\begin{eqnarray*}
&&
2 \sigma_0 \sum_{j=1}^p
\biggl(\sqrt{ 2( s_j + \tilde s_j ) \over n } + \sqrt{ 2( t+ \tilde s_j
+ 2 \log p ) \over n} \biggr) \bigl\| X
\bigl( \beta_j - \tilde\beta_j^0 \bigr)
\bigr\|_n
\\
&&\qquad
\le\delta_1 \sum_{j=1}^p \bigl\|
X \bigl( \beta_j - \tilde\beta_j^0 \bigr)
\bigr\|_n^2 + {4
\sigma_0^2 (s + \tilde s) \over n \delta_1 } + { 4\sigma_0^2
(t + 2 \log p )(s+ \tilde s) \over n \delta_1 },
\end{eqnarray*}
where $s = \sum_{j=1}^p s_j $, $\tilde s = \sum_{j=1}^p \tilde s_j $.
\end{pf}

\subsubsection{The set ${\cal T}_2$}\label{T2section}

\begin{theorem} \label{T2theorem} Assume Condition \ref
{edgescondition}. Then for all $t > 0 $,
\begin{eqnarray*}
&&
\PP\Biggl(\exists\pi\dvtx  \sum_{j=1}^p \biggl(
{ \| \tilde\epsilon_j
(\pi) \|_n^2 - | \tilde\omega_j^0 (\pi)|^2 \over| \tilde\omega_j^0
(\pi) |^2 } \biggr)^2
\\
&&\qquad\ge8 \biggl( { p t + (1+ 8 \tilde\alpha) \tilde s(\pi) \log p + 2 p
\log p \over n } \biggr) +8 \biggl( { 4 p (t^2 + \log^2 p) \over n^2
}
\biggr) \Biggr)
\\
&&\qquad\le2 \exp[-t].
\end{eqnarray*}
\end{theorem}

\begin{pf}
Define
\[
Z_j (\pi):= { \| \tilde\epsilon_j (\pi) \|_n^2 - | \tilde\omega_j^0
(\pi)|^2 \over| \tilde\omega_j^0 (\pi) |^2 }.
\]
Using the same argument as in (\ref{Bernstein}), we see
that for each $\pi$, and for all $t > 0$,
\[
\PP\biggl( \bigl| Z_j (\pi) \bigr| \ge2 \biggl( \sqrt{ t \over n}
+ {t \over
n} \biggr) \biggr) \le2 \exp[-t].
\]
Define $ {\mathbf Z}_j (\pi):=
{| Z_j (\pi) | / 2 a_j (\pi) } $,
where
\begin{eqnarray*}
a_j (\pi) &=& \biggl( \sqrt{ t + \tilde s_j (\pi) \log p + \log(1+ p) + \log p
\over n}
\\
&&\hspace*{5.2pt}{}+ {t + \tilde s_j (\pi) \log
p + \log(1+p) + \log p \over n} \biggr).
\end{eqnarray*}
It follows that
\begin{eqnarray*}
&&
\PP\Bigl( \max_{1 \le j \le p} \max_{0 \le m \le p } \max
_{\pi\in\Pi
_j (m) } {\mathbf Z}_j (\pi) \ge1 \Bigr)
\\
&&\qquad
\le2 p (p+1) \pmatrix{p \cr m} \exp\bigl[ - \bigl( t + m \log p + \log(1+p
) + \log
p \bigr) \bigr] \\
&&\qquad\le2 \exp[-t].
\end{eqnarray*}
Invoking $\log(1+ p) \le2 \log p $, we see that with probability at
least $1- 2\exp[-t]$, it holds that for all permutations $\pi$ and all $j$,
\[
\bigl| Z_j (\pi) \bigr| \le2 \sqrt{ t + \tilde s_j (\pi) \log p + 2 \log p
\over n} +
{t + \tilde s_j (\pi) \log p + 2 \log p \over n},
\]
which implies
\begin{eqnarray*}
\sum_{j=1}^p \bigl|Z_j (\pi)
\bigr|^2 &\le&4 \sum_{j=1}^p \biggl(
\sqrt{ t +
\tilde s_j (\pi) \log p +2 \log p \over n} + {t + \tilde s_j (\pi) \log
p + 2 \log p \over n} \biggr)^2
\\
&\le&8 \biggl( { p t + \tilde s(\pi) \log p + 2 p \log p \over n } \biggr
) \\
&&{}+8 \biggl( { 4 p t^2 + 8 \sum_{j=1}^p \tilde s_j^2 (\pi) \log
^2 p + 4 p \log^2 p \over n^2 }
\biggr).
\end{eqnarray*}
Next, we insert that for all $j$, $\tilde s_j (\pi) \le\tilde\alpha
n / (\log p ) $, to find
\[
\sum_{j=1}^p \tilde s_j^2
(\pi) \log^2 p \le\sum_{j=1}^n
\bigl( \tilde\alpha n/ (\log p ) \bigr) \tilde s_j (\pi)
\log^2 p = \tilde\alpha\tilde s (\pi) n \log p.
\]
We then arrive at
\begin{eqnarray*}
\sum_{j=1}^p \bigl|Z_j (\pi)
\bigr|^2 &\le&8 \biggl( { p t + (1+ 8 \tilde\alpha
) \tilde s(\pi) \log p + 2 p \log p \over n } \biggr)\\
&&{} +8 \biggl(
{ 4
p( t^2 + \log^2 p) \over n^2 } \biggr).
\end{eqnarray*}
\upqed\end{pf}

\subsubsection{The sets ${\cal T}_3 $ and ${\cal T}_0$}\label{T3T0section}

\begin{theorem} \label{normtheorem}
Assume Conditions~\ref{sigmacondition} and~\ref{Lambda-mincondition}.
For all $t>0$, with probability at least $1- 2 \exp[-t] $,
\[
\| X \beta\|_n \ge\biggl[ 3 \Lambda_{\mathrm{min}} /4 - \sqrt
{ 2( t+
\log p )\over n } - 3 \sigma_0 \sqrt{ s_{\beta} \log p \over n }
\biggr] \| \beta\|_2,
\]
uniformly in all $\beta\in\R^p$.
\end{theorem}

\begin{pf}
We follow here the arguments used in \citet{raskutti2010restricted},
which we slightly adjust to
the style of the present paper.
They show that for $\delta_3^{\prime}= 1/4$ [in fact for $\delta
_3^{\prime} = o(1)$ as $ n \rightarrow\infty$], and for all $r > 0$,
\[
\EE\inf_{\| \beta\|_1 \le r, \| X \beta\| = 1 } \| X \beta\|_n \ge1-
\delta_3^{\prime} - 3 \sigma_0 \sqrt
{ \log p \over n } r.
\]
Hence, for all $1 \le m \le p$,
\[
\EE\inf_{s_{\beta} \le m, \| \beta\|_2= 1 } \| X \beta\|_n \ge\bigl(1-
\delta_3^{\prime} \bigr) \Lambda_{\mathrm{min}} - 3
\sigma_0 \sqrt{ m \log p
\over n}.
\]
Apply the concentration inequality given in \citet
{massart1896concentration} to find that for all $t >0$,
\[
\PP\biggl( \Bigl[ \EE\inf_{s_{\beta} \le m, \| \beta\|_2 = 1
} \| X \beta\|_n
\Bigr] - \Bigl[ \inf_{s_{\beta} \le m, \|\beta\|_2 = 1 } \| X \beta\|_n
\Bigr] \ge
\sqrt{ 2t \over n } \biggr) \le2 \exp[-t].
\]
Thus
\[
\PP\biggl( \biggl[ \bigl(1- \delta_3^{\prime} \bigr)
\Lambda_{\mathrm{min}} - 3 \sigma_0 \sqrt{ m \log p \over n }
\biggr] - \Bigl[ \inf_{s_{\beta} \le m, \| \beta\|_2 = 1 } \| X \beta\|
_n \Bigr]
\ge\sqrt{ 2t \over n } \biggr) \le2 \exp[-t]
\]
and hence
\begin{eqnarray*}
&&
\PP\biggl( \exists\beta\dvtx  \biggl[ \bigl(1- \delta_3^{\prime} \bigr)
\Lambda_{\mathrm{min}} - 3 \sigma_0 \sqrt{ s_{\beta} \log p \over n }
\biggr] \| \beta\|_2 - \| X \beta\|_n
\ge\sqrt{ 2( t+ \log p )\over n } \| \beta\|_2 \biggr)\qquad \\
&&\qquad\le2 \exp[-t].
\end{eqnarray*}
\upqed\end{pf}

\begin{lemma} \label{omegalemma}
Assume Conditions~\ref{sigmacondition},~\ref{Lambda-mincondition},
\ref{hatedgescondition} and~\ref{edgescondition} and that
\[
1/ K_0:= 3 \Lambda_{\mathrm{min}} /4 - \sqrt{ 2( t+ \log p )\over n }
- 3 \sigma_0 \sqrt{ \alpha+ \tilde\alpha} >0.
\]
Let for some $t > 0$,
\[
\tilde{\cal T}_3:= \biggl\{ \| X \beta\|_n \le\biggl[ 3
\Lambda_{\mathrm{min}} /4 - \sqrt{ 2( t+ \log p )\over n } - 3
\sigma_0 \sqrt{
s_{\beta} \log p \over n } \biggr] \| \beta
\|_2, \forall\beta\biggr\}.
\]
Then $\PP( \tilde{\cal T}_3 ) \ge1- 2 \exp[-t]$ and one has on
${\cal T}_3$, for all $B = ( \beta_1,\ldots, \beta_p) \in{\cal B} $
and all $\pi$ and all $j$,
%
\begin{equation}
\label{eigenvalue} \bigl\| X \bigl( \beta_j - \tilde\beta_j^0
(\pi) \bigr) \bigr\|_n \ge\bigl\| \beta_j - \tilde
\beta_j^0 (\pi) \bigr\|_2 / K_0^2.
\end{equation}
Moreover, on $\bar{\cal T}_3$, also $ \hat\omega_j^2 \ge1/ K_0 ^2 $
for all $j$.
\end{lemma}

\begin{pf}
Theorem~\ref{normtheorem} states that $\PP( \tilde
{\cal T}_3 ) \ge1- 2 \exp[-t]$.
Result (\ref{eigenvalue}) follows immediately, since $s_{\beta_j} +
s_{\tilde\beta_j^0} \le( \alpha+ \tilde\alpha) n / \log p $.
For the last result, we
define
$\hat\beta_{k,j}^{-}:= - \hat\beta_{k,j}$ for
$k \not= j$ and $\hat\beta_{j,j}^{-} = 1$. Then on $\tilde{\cal T}_3 $,
\[
\hat\omega_j^2 = \bigl\| X \hat\beta_{j}^{-}
\bigr\|_n^2 \ge\bigl\| \hat\beta_{j}^{-}
\bigr\|_2^2 / K_0^2 \ge1/ K_0^2.
\]
\upqed\end{pf}

\subsection{Collecting the results}\label{collectingsection}
%
\begin{lemma} \label{exactlemma}
Define $( \tilde B_0, \tilde\Omega_0 ):=
( B_0 ( \hat\pi), \Omega_0 ( \hat\pi) ) $.
Assume
Conditions~\ref{sigmacondition},\break \ref{Lambda-mincondition},
\ref{hatedgescondition},
\ref{edgescondition} and~\ref{beta-mincondition}.
Suppose we are on $\bigcap_{k=0}^3 {\cal T}_k $ with
$0 < \delta_1 < 1/ K_0^2 $ and $0< \delta_2 < 1/ (2K_0^4 \sigma_0^4 ) $ and
$\delta_3 - \lambda_3 \sqrt{ \alpha+ \tilde\alpha} \sqrt{ n / \log
p } \ge1/ K_0 > 0 $.
Take the tuning parameter $\lambda^2 > \lambda_1^2 / \delta_1+ \lambda
_2^2 /\delta_2 $. Let
\begin{eqnarray*}
\delta_B&\le&{1 \over K_0^2} \biggl({1 \over K_0^2}
- \delta_1 \biggr),\qquad \delta_W\le{1 \over K_0^2 }
\biggl( { 1 \over2 K_0^4 \sigma_0^4 }- \delta_2 \biggr),
\\
\delta_s &\le&\biggl( 1 - { \lambda_1^2 \over\lambda^2 \delta_1 }-
{
\lambda_2^2 \over\lambda^2 \delta_2 } \biggr),\qquad \lambda_0^2:=
{ ({p / s_0} + 1 ) \lambda_2^2 \over\delta_2 } + {
\lambda_1^2 \over\delta_1 }.
\end{eqnarray*}

Let $\tilde\lambda^2 \delta_B: = \lambda^2+ \lambda_0^2 $, and
$\eta_2^2: = \eta_0^2 \tilde\lambda^2 n/ \log p = \eta_0^2 ( \lambda
^2 + \lambda_0^2 ) ( n/ \log p )/ \delta_B^2 $.
Assume
\[
\biggl( \lambda^2 \delta_s - { \lambda_0^2 \over1- \eta_1 - \eta_2^2}
\biggr):= \lambda^2 \delta_{\eta} >0.
\]
Then
\[
\delta_B \|\hat B - \tilde B_0 \|_{F}^2
+ \delta_W \| \hat\Omega- \tilde\Omega_0
\|_F^2 + \lambda^2 \delta_{ \eta} \hat
s \le\lambda^2 s_0
\]
and $\hat s \ge(1- \eta_1 - \eta_2^2 ) \tilde s \ge(1- \eta_1 - \eta
_2^2 ) s_0 $.
\end{lemma}

\begin{pf}
This follows from Lemmas~\ref{startlemma} and
\ref{beta-min2lemma}.
\end{pf}

\begin{lemma}\label{lambdalemma} Assume
Conditions~\ref{sigmacondition},~\ref{Lambda-mincondition},
\ref{hatedgescondition} and~\ref{edgescondition}, with
\[
3 \Lambda_{\mathrm{min}}/ 4 - \sqrt{ 2(t+ \log p ) \over n} - 3
\sigma_0 \sqrt{\alpha+ \tilde\alpha}\ge1/ K_0 >0.
\]
Take
\begin{eqnarray*}
\lambda_1^2 &:=& { 4 \sigma_0^2 (1+ t + 2 \log p ) \over\log p }
{ \log
p \over n},
\\
\lambda_2^2 &:=& 4 \biggl( 3+ 8 \tilde\alpha+
{ (t+ 2 \log p ) \over
\log p } + { 4 ( t^2 + \log^2 p ) \over n \log p } \biggr) {\log p
\over n},
\\
\lambda_3^2 &:=& 9 \sigma_0^2
{\log p \over n},\qquad \delta_3:= {3 \over4}
\Lambda_{\mathrm{min}} - \sqrt{ 2(t+ \log p ) \over
n }.
\end{eqnarray*}
Then
\[
\PP\Biggl( \bigcap_{k=0}^3 {\cal
T}_k \Biggr) \ge1- 5 \exp[-t].
\]
\end{lemma}

\begin{pf}
This follows from combining Theorems~\ref{T1theorem},
\ref{T2theorem} and
Lemma~\ref{omegalemma}.
\end{pf}

\begin{theorem}\label{exact2theorem}
Assume
Conditions~\ref{sigmacondition},~\ref{Lambda-mincondition},
\ref{hatedgescondition},~\ref{edgescondition} and~\ref{beta-mincondition}.
Let us take $t = \log p$, giving $\alpha_0 = 4/p $ (suppose $p$ is large).
Take $n$ sufficiently large, and $\log p / n$ bounded.
Let
\begin{eqnarray*}
c_1:\!&=& 96,\qquad c_2:= 3840,\\
c&=& 4 \biggl(
{ ( p / s_0 +1 ) c_2 \sigma_0^2 \over\Lambda_{\mathrm{min}}^4 } + { c_1
\sigma_0^2 \over\Lambda_{\mathrm{min}}^2 } \biggr) + 2 \biggl(
{ c_1 \sigma_0^2 \over\Lambda_{\mathrm{min}}^2 } + {c_2 \sigma
_0^2 \over\Lambda_{\mathrm{min}}^4 } \biggr).
\end{eqnarray*}
Some possible choices for the constants are
\begin{eqnarray*}
\alpha &=& \tilde\alpha= { \Lambda_{\mathrm{min}}^2 \over288 \sigma_0^2 },\qquad K_0:=
{ 2 \over\Lambda_{\mathrm{min}}},
\\
\delta_1 :\!&=& {\Lambda_{\mathrm{min}}^2\over8 },\qquad \delta_2:=
{ \Lambda_{\mathrm{min}}^4 \over64 \sigma_0^4 },\qquad \delta_3 := { \Lambda_{\mathrm{min}} \over2
},
\\
\lambda^2:\!&=& c { \log p \over n },\qquad \lambda_1^2
= 12 \sigma_0^2{ \log
p \over n },\qquad
\lambda_2^2 = 60 {\log p \over n},\qquad
\lambda_3^2 = 9 \sigma_0^2
{\log p \over n}.
\end{eqnarray*}

Then
\[
\delta_B = { \Lambda_{\mathrm{min}}^4 \over32},\qquad \delta_W=
{ \Lambda_{\mathrm{min}}^6 \over256 \sigma_0^4 },\qquad \delta_s = \biggl( 1- {
c_1 \sigma_0^2 \over c \Lambda_{\mathrm{min}}^2 }
- {c_2\sigma_0^4 \over c \Lambda_{\mathrm{min}}^4 } \biggr).
\]
We let
\begin{eqnarray*}
\lambda_0^2 :\!&=& \biggl( { ( p / s_0 +1 ) c_2 \sigma_0^4 \over\Lambda
_{\mathrm{min}}^4 } +
{ c_1 \sigma_0^2 \over\Lambda_{\mathrm{min}}^2 } \biggr) {\log p \over n
},
\\
\tilde\lambda^2 &=& \biggl( c + { ( p / s_0 +1 ) c_2 \sigma_0^4 \over
\Lambda_{\mathrm{min}}^4 } +
{ c_1 \sigma_0^2 \over\Lambda_{\mathrm{min}}^2 } \biggr) {32 \over\Lambda
_{\mathrm{min}}^4 } {\log p \over n}
\end{eqnarray*}
and
\begin{eqnarray*}
\delta_{\eta} &=& {1 \over2},\qquad \eta_1 =0,\qquad 2
\eta_0^2= \biggl( c + { ( p / s_0 +1 ) c_2 \sigma_0^4 \over\Lambda_{\mathrm{min}}^4 } +
{ c_1
\sigma_0^2 \over\Lambda_{\mathrm{min}}^2 } \biggr)^{-1},
\\
\eta_2^2 &=& \eta_0^2 \biggl( c +
{ ( p / s_0 +1 ) c_2 \sigma_0^4 \over
\Lambda_{\mathrm{min}}^4 } + { c_1 \sigma_0^2 \over\Lambda_{\mathrm{min}}^2 }
\biggr) {32 \over\Lambda_{\mathrm{min}}^4 }.
\end{eqnarray*}
\end{theorem}

\begin{pf}
This follows from using some bounds and exact choices in Lemmas
\ref{exactlemma} and~\ref{lambdalemma}.
In particular, with $\lambda^2 = c \log p / n $, we take $\tilde
\lambda^2 = (\lambda^2 + \lambda_0^2 )/ \delta_B$.
With $\eta_1=0$ and $\delta_{\eta} = 1/2$, the equation
\[
\biggl( \lambda^2 \delta_s - {\lambda_0^2 \over1- \eta_2^2 }
\biggr) = {\lambda^2 \over2 }
\]
gives
\begin{eqnarray*}
&&
c \biggl( 1- { c_1 \sigma_0^2 \over c \Lambda_{\mathrm{min}}^2 } - {c_2
\sigma_0^4 \over c \Lambda_{\mathrm{min}}^4 } \biggr)
\\
&&\qquad{}- \biggl( { ( p / s_0 +1 ) c_2 \sigma_0^4 \over\Lambda_{\mathrm{min}}^4 } +
{ c_1 \sigma_0^2 \over\Lambda_{\mathrm{min}}^2 } \biggr)\\
&&\qquad\hspace*{11.4pt}{}\times \biggl(1-
\eta_1 - \eta_0^2 \biggl( c +
{ ( p / s_0 +1 ) c_2 \sigma
_0^4 \over\Lambda_{\mathrm{min}}^4 } + {c_1 \sigma_0^2 \over\Lambda_{\mathrm{min}}^2 } \biggr)
\biggr)^{-1}
= {c \over2}.
\end{eqnarray*}
With
\[
2 \eta_0^2= \biggl( c + { c_2 \sigma_0^4 \over\Lambda_{\mathrm{min}}^4 } +
{ c_1 \sigma_0^2 \over\Lambda_{\mathrm{min}}^2 } \biggr)^{-1}
\]
we have to solve for $c$
\[
c \biggl( 1- { c_1 \sigma_0^2 \over c \Lambda_{\mathrm{min}}^2 } - {c_2
\sigma_0^4 \over c \Lambda_{\mathrm{min}}^4 } \biggr) - 2
\biggl( { ( p / s_0 +1 ) c_2 \sigma_0^4 \over\Lambda_{\mathrm{min}}^4 } + {
c_1 \sigma_0^2 \over\Lambda_{\mathrm{min}}^2 } \biggr) =
{ c \over2}.
\]
This yields
\[
c= 4 \biggl( { ( p / s_0+1 ) c_2 \sigma_0^4 \over\Lambda_{\mathrm{min}}^4
} + { c_1 \sigma_0^2 \over\Lambda_{\mathrm{min}}^2 } \biggr) + 2
\biggl( { c_1 \sigma_0^2 \over\Lambda_{\mathrm{min}}^2 } + {c_2 \sigma
_0^4 \over\Lambda_{\mathrm{min}}^4 } \biggr).
\]
\upqed\end{pf}

\subsection{\texorpdfstring{Proof of Theorem \protect\ref{equalvariancestheorem}}{Proof of Theorem 5.1}}
\label{equalvariancesproofsection}

We investigate what happens on the set $\bigcap_{k=1}^3 {\cal T}_3$
defined in Section
\ref{subsetprobabilitysection}. The results in Section \ref
{T-section} say that
$\bigcap_{k=1}^3{\cal T}_k$ has probability at least $4 \exp[-t]$ for a proper
choice of the constants and parameters involved. Theorem \ref
{equalvariancestheorem}
then follows directly.

\begin{lemma} \label{equalvarianceslemma}
Assume Conditions~\ref{sigmacondition},~\ref{omega-mincondition} and
\ref{psmallcondition}. Suppose we are on $\bigcap_{k=1}^3 {\cal T}_k$,
with
\begin{eqnarray*}
\lambda^2 &>& \lambda_1^2 /
\delta_1,
\\
\delta_3 - \lambda_3 \sqrt{ 2 \alpha_* } &\ge&
{1 \over K_0 } >0
\end{eqnarray*}
and
%
\begin{equation}
\label{omega-minequation} \biggl( {1 \over2 \sigma_0^4 } -
\delta_2 \biggr) \ge\eta_{\omega
} \biggl( { 2 \lambda_1^2 \over\delta_2}
+ {\lambda_1^2 \over\delta_1} + \lambda^2 \biggr) {\alpha_* n \over
\log p }.
\end{equation}
Then $\hat\pi= \pi_0$ and
\[
\biggl( { 1- \delta_1 \over K_0^2 } \biggr) \| \hat B - B_0
\|_F^2 + \biggl( \lambda^2 -
{\lambda_1^2 \over\delta_1 } \biggr) \hat s \le\biggl( \lambda^2 +
{ \lambda_1^2 \over
\delta_1 } \biggr) s_0.
\]
\end{lemma}

\begin{pf}
We have
\[
\sum_{j=1}^p \| X_j - X \hat
\beta_j \|_n^2 + \lambda^2 \hat s
\le\sum_{j=1}^p \| \epsilon_j
\|_n^2 + \lambda^2 s_0 = \sum
_{j=1}^p{ \|
\tilde\epsilon_j \|_n^2 \over
| \tilde\omega_j^0 |^2 }.
\]
So we find
\begin{eqnarray*}
&&
\sum_{j=1}^p \bigl\| X \bigl( \hat
\beta_j - \tilde\beta_j^0\bigr)
\bigr\|_n^2 + \lambda^2 \hat s \\
&&\qquad\le 2 \sum
_{j=1}^p \tilde\epsilon_j^T
X \bigl( \hat\beta_j - \tilde\beta_j^0
\bigr) / n + \sum_{j=1}^p \| \tilde
\epsilon_j \|_n^2 \biggl(
{1 \over| \tilde
\omega_j^0 |^2 } - 1 \biggr) + \lambda^2 s_0.
\end{eqnarray*}

We have
\[
\sum_{j=1}^p \| \tilde
\epsilon_j \|_n^2 \biggl(
{1 \over| \tilde
\omega_j^0 |^2 } - 1 \biggr)= \sum_{j=1}^p
{ \| \tilde\epsilon_j \|_n^2 - |\tilde\omega_j^0 |^2
\over
| \tilde\omega_j^0 |^2 } \bigl( \bigl| \tilde\omega_j^0
\bigr|^2 - 1 \bigr) + \sum_{j=1}^p
\bigl( 1- \bigl|\tilde\omega_j^0 \bigr|^2 \bigr).
\]

We know that
\[
\log\bigl( \operatorname{det} ( \Sigma_0 )\bigr) = \sum
_{j=1}^p \log\bigl|\tilde\omega_j^0
\bigr|^2 = \sum_{j=1}^p \log\bigl|
\omega_j^0 \bigr|^2 = 0,
\]
since $|\omega_j^0 |^2 = 1$ for all $j$. Moreover
\[
\log(1+x) \le x - {1 \over2 (1+c)^2 } x^2,\qquad -1 < x \le c.
\]
So, since $| \tilde\omega_j^0 |^2 \le\sigma_0^2 $,
\[
\log\bigl| \tilde\omega_j^0\bigr|^2 \le\bigl( {\bigl|
\tilde\omega_j^0 \bigr|^2} - 1 \bigr) -
{1
\over2 \Sigma_0^4} \bigl({\bigl|\tilde\omega_j^0
\bigr|^2 } - 1 \bigr)^2.
\]
Hence
\[
0 \le\sum_{j=1}^p \bigl( \bigl|\tilde
\omega_j^0 \bigr|^2 - 1 \bigr) -
{1 \over2 \sigma
_0^4} \sum_{j=1}^p
\bigl( \bigl|\tilde\omega_j^0 \bigr|^2 - 1
\bigr)^2.
\]
This gives
\[
\sum_{j=1}^p \bigl(1- \bigl|\tilde
\omega_j^0 \bigr|^2 \bigr)\le-
{1 \over2 \sigma_0^4} \sum_{j=1}^p
\bigl( \bigl|\tilde\omega_j^0 \bigr|^2 - 1
\bigr)^2.
\]
Therefore
\begin{eqnarray*}
&&
\sum_{j=1}^p \bigl\| X \bigl( \hat
\beta_j - \tilde\beta_j^0\bigr)
\bigr\|_n^2 +{1 \over2
\sigma_0^4} \sum
_{j=1}^p \bigl( \bigl|\tilde\omega_j^0
\bigr|^2 - 1 \bigr)^2 + \lambda^2 \hat s
\\
&&\qquad
\le2 \sum_{j=1}^p \tilde
\epsilon_j^T X \bigl( \hat\beta_j - \tilde
\beta_j^0 \bigr) / n + \lambda^2
s_0 +\sum_{j=1}^p
{ \| \tilde\epsilon_j \|_n^2 -
|\tilde\omega_j^0 |^2 \over
| \tilde\omega_j^0 |^2 } \bigl( \bigl| \tilde\omega_j^0
\bigr|^2 - 1 \bigr)
\\
&&\qquad\le2 \sum_{j=1}^p \tilde
\epsilon_j^T X \bigl( \hat\beta_j - \tilde
\beta_j^0 \bigr) / n + \lambda^2
s_0 + { \lambda_2^2 p \over\delta_2 } + \delta_2 \sum
_{j=1}^p \bigl( \bigl|\tilde\omega_j^0
\bigr|^2 - 1 \bigr)^2,
\end{eqnarray*}
where we invoked that we are on the set ${\cal T}_2$.
We find
\begin{eqnarray*}
&&
\sum_{j=1}^p \bigl\| X \bigl( \hat
\beta_j - \tilde\beta_j^0\bigr)
\bigr\|_n^2 + \biggl( {1
\over2 \sigma_0^4} -
\delta_2 \biggr) \sum_{j=1}^p
\bigl( \bigl|\tilde\omega_j^0 \bigr|^2 - 1
\bigr)^2 + \lambda^2 \hat s
\\
&&\qquad
\le2 \sum_{j=1}^p \tilde
\epsilon_j^T X \bigl( \hat\beta_j - \tilde
\beta_j^0 \bigr) / n + \lambda^2
s_0 +{\lambda_2^2 p \over\delta_2 }.
\end{eqnarray*}

This gives in a next step, using that we are on ${\cal T}_1$,
\begin{eqnarray*}
&&
( 1 - \delta_1 ) \sum_{j=1}^p
\bigl\| X \bigl( \hat\beta_j - \tilde\beta_j^0\bigr) \bigr\|
_n^2 + \biggl( {1 \over2 \sigma_0^4} -
\delta_2 \biggr) \sum_{j=1}^p
\bigl( \bigl|\tilde\omega_j^0 \bigr|^2 - 1
\bigr)^2 + \lambda^2 \hat s
\\
&&\qquad
\le{\lambda_2^2 ( p + \tilde s ) \over
\delta_2 } + { \lambda_1^2 ( \tilde s + \hat s) \over\delta_1} + \lambda^2
s_0.
\end{eqnarray*}
Hence, using that we are on ${\cal T}_3$ and invoking Condition \ref
{psmallcondition},
\begin{eqnarray*}
&&{1 \over K_0^2 }( 1 - \delta_1 ) \| \hat B - B_0
\|_F^2 + \biggl( {1
\over2 \sigma_0^4} -
\delta_2 \biggr) \sum_{j=1}^p
\bigl( \bigl|\tilde\omega_j^0 \bigr|^2 - 1
\bigr)^2 + \biggl( \lambda^2 - {\lambda_1^2 \over\delta_1 }
\biggr) \hat s
\\
&&\qquad
\le{\lambda_2^2 ( p + \tilde s ) \over
\delta_2 } + { \lambda_1^2 \tilde s \over\delta_1} + \lambda^2
s_0
\\
&&\qquad\le\biggl( { 2 \lambda_2^2 \over\delta_2} + {\lambda_1^2 \over\delta
_1 } +
\lambda^2 \biggr) p^2,
\end{eqnarray*}
where we use that $\tilde s \le p^2 $ and $s_0 \le p^2 $ (and also $p
\le p^2 $).
Since (using again Condition~\ref{psmallcondition}) $p \log p / n \le
\alpha_*$, and
\[
\sum_{j=1}^p \bigl( | \tilde
\omega_j |^2 -1 \bigr)^2 > p /\eta
_{\omega} \qquad\mbox{if } \hat\pi\not= \pi_0,
\]
find that if $\hat\pi\not= \pi_0$,
\[
\biggl( {1 \over2 \sigma_0^4} - \delta_2 \biggr)
{ p \over\eta
_{\omega} } < \biggl( { 2 \lambda_2^2 \over\delta_2} +
{\lambda_1^2
\over\delta_1 } + \lambda^2 \biggr) \alpha_*
{ n \over\log p },
\]
which is in contradiction with Condition~\ref{omega-mincondition} and
the further, condition (\ref{omega-minequation}) imposed in this lemma.
So we must have $\hat\pi= \pi_0$, and thus $\tilde\omega_j^0 = 1 $
for all~$j$.
The result now follows from restarting the proof with
$\tilde\omega_j^0 = 1$ for all $j$ plugged in.
\end{pf}




\printaddresses

\end{document}